\documentclass[12pt]{amsart}
\usepackage{epsfig}
\usepackage{amsmath}
\usepackage{amsfonts}
\usepackage{amssymb}
\usepackage{mathabx}
\usepackage{epigraph, fancyhdr}

\def\C{\mathcal C}
\def\D{\mathcal D}
\def\F{\mathcal F}
\def\K{\mathcal K}

\def\H{\mathcal H}
\def\O{\mathcal O}
\def\T{\mathcal T}

\def\Z{\mathcal Z}
\def\1{\mathbf 1}
\def\M{\mathcal M}

\def\QQ{\mathbb Q}
\def\ZZ{\mathbb Z}
\def\CC{\mathbb C}

\def\Res{\operatorname{Res}}

\def\hat{\widehat}
\def\tilde{\widetilde}
\def\p{\partial}
\def\a{\alpha}
\def\ga{\alpha}
\def\b{\beta}

\def\f{{\mathbf f}}
\def\g{{\mathbf g}}
\def\t{{\mathbf t}}
\def\q{{\mathbf q}}
\def\pp{{\mathbf p}}
\def\x{{\mathbf x}}

\def\ll{{\mathbf l}}
\def\vv{{\mathbf v}}
\def\uu{{\mathbf u}}
\def\ww{{\mathbf w}}

\def\gs{\sigma}

\def\la{\lambda}

\def\gL{\Lambda}
\def\gG{\Gamma}

\def\gS{\Sigma}

\def\m{{\mathfrak m}}
\def\h{\hbar}
\def\w{\wedge}
\def\lan{\langle}
\def\ran{\rangle}
\def\str{\operatorname{str}}

\def\ev{\operatorname{ev}}
\def\ft{\operatorname{ft}}

\def\td{\operatorname{td}}
\def\ch{\operatorname{ch}}
\def\qch{\operatorname{qch}}
\def\eu{\operatorname{eu}}

\def\Eu{\operatorname{Eu}}

\def\tr{\operatorname{tr}}
\def\str{\operatorname{str}}

\def\square{\Box}
\def\und{\underline}
\renewcommand{\Delta}{\triangle}
\makeatletter
\DeclareFontFamily{OMX}{MnSymbolE}{}
\DeclareSymbolFont{MnLargeSymbols}{OMX}{MnSymbolE}{m}{n}
\SetSymbolFont{MnLargeSymbols}{bold}{OMX}{MnSymbolE}{b}{n}
\DeclareFontShape{OMX}{MnSymbolE}{m}{n}{
    <-6>  MnSymbolE5
   <6-7>  MnSymbolE6
   <7-8>  MnSymbolE7
   <8-9>  MnSymbolE8
   <9-10> MnSymbolE9
  <10-12> MnSymbolE10
  <12->   MnSymbolE12
}{}
\DeclareFontShape{OMX}{MnSymbolE}{b}{n}{
    <-6>  MnSymbolE-Bold5
   <6-7>  MnSymbolE-Bold6
   <7-8>  MnSymbolE-Bold7
   <8-9>  MnSymbolE-Bold8
   <9-10> MnSymbolE-Bold9
  <10-12> MnSymbolE-Bold10
  <12->   MnSymbolE-Bold12
}{}

\let\llangle\@undefined
\let\rrangle\@undefined
\DeclareMathDelimiter{\llan}{\mathopen}%
                     {MnLargeSymbols}{'164}{MnLargeSymbols}{'164}
\DeclareMathDelimiter{\rran}{\mathclose}%
                     {MnLargeSymbols}{'171}{MnLargeSymbols}{'171}
\makeatother

\title[Quantum HRR in all genera]
      {Permutation-equivariant \\ quantum K-theory IX. \\
      Quantum Hirzebruch-Riemann-Roch \\ in all genera}

\author[A. Givental]{Alexander GIVENTAL}
\thanks{This material is based upon work supported by the National 
Science Foundation under Grant  DMS-1611839, and by the IBS Center for Geometry 
and Physics, POSTECH, Korea.}

\date{August 24, 2017}

\begin{document}

\begin{abstract}

We introduce the most general to date version of the permutation-equivariant quantum K-theory, and express its total descendant potential in terms of cohomological Gromov-Witten invariants. This is the higher-genus analogue of adelic characterization \cite{GiTo}, and is based on the application of the Kawasaki-Riemann-Roch formula \cite{Kaw} to moduli spaces of stable maps.    
  
\end{abstract}

\maketitle

\section*{Introduction}

Cohomological Gromov-Witten invariants of a compact K\"ahler manifold $X$
are defined as various intersection numbers in moduli spaces of stable maps,
denoted here $X_{g,n,d}$ with $g$, $n$, $d$ standing for the genus, number of marked points, and degree of the maps. The K-theoretic counterpart of GW-theory studies holomorphic Euler characteristics of appropriate vector bundles over the moduli spaces. The action of permutations of the marked points on the sheaf cohomology of such bundles leads to the refined version of
the theory, which we call permutation-equivariant. In genus $0$, a complete description of K-theoretic GW-invariants in terms of cohomological ones was obtained in \cite{GiTo}, and then applied to the permutation-equivariant
theory in the previous papers of the present series (see Part III or Part VII). 

Conceptually the cohomological description of K-theoretic invariants is based on Kawasaki's version of Hirzebruch--Riemann--Roch formula \cite{Kaw} (or more precisely, its virtual variant \cite{ToK}) applied on the moduli spaces $X_{g,n,d}$. An early version of this approach to the higher genus problem is used in the preprint \cite{ToH} by V. Tonita. I am thankful to him for numerous discussions and corrections. 

As it was found in \cite{GiTo}, in genus $0$ the solution can be described in the form of {\em adelic characterization}. Roughly speaking, genus-$0$ K-theoretic GW-invariants of $X$ are encoded by a certain Lagrangian cone in a symplectic space whose elements are rational functions in one complex variable, $q$, with vector values in $K^0(X)$. The adelic characterization says that a rational function lies in the cone if and only if the Laurent series expansion of it at each root of unity $q=\zeta$ passes a certain test. Namely, the expansion (as an element in the symplectic space of Laurent series with coefficients in $K^0(X)$) should represent certain cohomological GW-invariants of the {\em orbifold} target space $X/\ZZ_M$, where $M$ is the order of $\zeta$ as a root of unity.

This paper establishes the higher genus version of adelic characterization. It involves quantization of the aforementioned symplectic formalism. In this Introduction, we don't give a complete formulation of the ultimate theorem (because it requires so many poorly motivated ingredients and notations, that the resulting formula, we fear, would become incomprehensible), but merely outline the quantum-mechanical structure of the adelic formula relating K-theoretic GW-invariants with cohomological ones.

A thorough definition of the permutation-equivariant GW-invariants and of the appropriate generating functions will be given in Section 1. In Section 2, we sketch the geometric machinery which shows, in principle, how to reduce the computation of K-theoretic to cohomological GW-invariants. In Sections 3 and 4, we describe the language of symplectic loop spaces and their Fock spaces where various generating functions for GW-invariants live. Using this language, we will accurately build the ingredients of the ultimate formula starting from cohomological GW-invariants. The remaining details of the proof will be provided in Sections 5--9.    

\medskip

By definition, permutation-equivariant K-theoretic GW-invariants take values in
a ground coefficient ring, $\Lambda$, which is a $\lambda$-algebra, i.e. is equipped with the action of Adams operations $\Psi^r:\Lambda \to \Lambda$, $r=1,2,3,\dots$, which are ring homomorphisms from $\Lambda$ to itself, and satisfy $\Psi^1=\operatorname{id}$, $\Psi^r\Psi^s=\Psi^{rs}$.

The {\em total descendant potential} $\D_X$ for permutation-equivariant GW-invariants of $X$ is defined (in Section 1) as a $\Lambda$-valued function of a sequence $\t=(\t_1,\t_2,\dots,\t_r,\dots) $ of Laurent polynomials\footnote{Foreshadowing the definition let us mention here that $\t_r$ will be used as the input in the correlators of permutation-equivariant quantum K-theory at those marked points which belong to cycles of length $r$ in the cycle structure of the permutation.}  in $q$ with vector coefficients in $K^0(X)\otimes \Lambda$. 
It also depends on the ``Planck constant'' $\h$, and can be interpreted as
an element of the Fock space associated with a certain symplectic space
$(\K^{\infty}, \Omega^{\infty})$. 
  
Namely, put $K:=K^0(X)\otimes \Lambda$, and consider the space $\K$ or rational $K$-valued functions of $q$ which are allowed to have poles only at $q=0,\infty$, or at roots of unity. Equip $\K$ with the $\Lambda$-valued symplectic form
\[ \Omega (\f,\g):= -\left[\Res_{q=0}+\Res_{q=\infty}\right] (\f(q^{-1}), \g(q))\, \frac{dq}{q},\]
where $(a,b):=\chi (X; a\otimes b)$ is the K-theoretic Poincar\'e pairing
on $K$, and with the Lagrangian polarization $\K=\K_{+}\oplus \K_{-}$, where
\[ \K_{+}:=K[q,q^{-1}], \ \ \K_{-}:=\left\{ \f \in \K \, | \, \f(\infty) = 0,\  \f(0)\neq \infty \right\} .\]
By definition, $\K^{\infty}$ consists of sequences $\f=(\f_1,\f_2,\dots,\f_r,\dots)$ of elements of $\K$. It is equipped with the symplectic form 
\[ \Omega^{\infty}(\f,\g):=\sum_{r=1}^{\infty} \frac{\Psi^r}{r}\, \Omega (\f_r, \g_r),\]
and Lagrangian polarization $\K_{\pm}^{\infty}=\{ \f=(\f_1,\f_2,\dots)\ |\ \forall r, \f_r\in \K_{\pm}\}$.
The total descendant potential $\D_X$, which is naturally a function of $\t=(\t_1,\t_2,\dots) \in \K_{+}^{\infty}$ (depending on the parameter $\h$), is considered as a function on $\K^{\infty}$ constant in the direction of $\K^{\infty}_{-}$, and in this capacity is interpreted as a ``quantum state'',
$\lan \D_X \ran$, an element of the Fock space associated with $(\K^{\infty},\Omega^{\infty})$.

On the cohomological side, for each $M=1,2,3,\dots$, let $\ZZ_M=\ZZ/M\ZZ$ denote the cyclic group of order $M$, and $\CC^{M-1}=\CC [\ZZ_M]/\CC$ be the quotient of the regular representation of $\ZZ_M$ by the trivial one. Over the global quotient orbifold $X/\ZZ_M$ (where the action of $\ZZ_M$ is trivial), introduce the orbibundle $T_X\otimes \CC^{M-1}$,
and denote by $E_M$ its total (orbi)space. What we need is a certain twisted cohomological GW-theory of $X/\ZZ_M$, which can be interpreted as the fake quantum K-theory\footnote{In fake quantum K-theory, genuine holomorphic Euler characteristics of orbibundles over moduli spaces of stable maps are replaced with their fake versions: $\chi^{fake}(\M; V):= \int_{\M} \ch (V) \td (T_{\M})$, and are therefore cohomological in nature.} of the non-compact orbifold $E_M$. Denote by $\D_{X/\ZZ_M}^{tw}$ the total descendant potential of such a theory.
Using a series of ``quantum Riemann-Roch theorems'' available in the literature
(see \cite{Co, CGL, JK, ToT, ToT, ToTs, Ts}), it will be shown in Sections 6,7 how to link this generating function directly to the total descendant potential $\D_X^H$ of the ordinary cohomological GW-theory of $X$. So we will assume here that all the functions $\D_{X/\ZZ_M}^{tw}$ are given.

Each $\D_{X/\ZZ_M}^{tw}$ can be considered as a quantum state, $\lan\D_{X/\ZZ_M}^{tw}\ran$, an element of the Fock space associated with the appropriate symplectic space, $(\K_{(M)}^{tw}, \Omega_{(M)}^{tw})$. This space is a direct sum of $M$ {\em sectors} corresponding to
$M$th roots of unity $\zeta$. Each sector is represented by the space $\K^{(\zeta)}$ isomorphic to the space $K((q-1))$ of vector-valued Laurent series in $q-1$. The symplectic form $\Omega_{(M)}^{tw}$ pairs $\K^{(\zeta)}$
with $\K^{(\zeta^{-1})}$ by the non-degenerate pairing
\[ (f,g) \mapsto \frac{1}{M}\Res_{q=1} (f(q^{-1}),g(q))^{(r)} \frac{dq}{q}.\]
It is based on the twisted Poincar\'e pairing on $K$ characterized by
\[ (\Psi^r a, \Psi^r b)^{(r)} = r \Psi^r(a,b),\]
where $r=r(\zeta)$ equals the index of the subgroup generated by $\zeta$ in the multiplicative group of all $M$th roots of unity. 

Note that when $M$ runs all positive integers, each root of unity $\zeta$ of {\em primitive} order $m=m(\zeta)$ occurs among $M$th roots of unity infinitely many times distinguished by the values of the index $r(\zeta)=M/m(\zeta)=1,2,3\dots$. Consequently the direct sum $\oplus_{M=1}^{\infty} \K_{(M)}^{tw}$ can be rearranged according to the indices $r$ into the {\em adelic space}
\[ \und{\K}^{\infty}:=\oplus_{\text{roots of unity $\zeta$}} \oplus_{r=1}^{\infty} \K_r^{(\zeta)}\]
(here $\K_r^{(\zeta)}$ is the $r$th copy of $\K^{(\zeta)}$) with the symplectic form
\[ \und{\Omega}^{\infty}(\f,\g)=\sum_{\zeta} \frac{1}{m(\zeta)}\sum_{r=1}^{\infty}\frac{1}{r} \Res_{q=1}(f_r^{(\zeta)}(q^{-1}),g_r^{(\zeta^{-1})}(q))^{(r)} \,  \frac{dq}{q}  .\]
Thus, the {\em adelic tensor product}
\[ \und{\D}_X:= \otimes_{M=1}^{\infty} \D_{X/\ZZ_M}^{tw}\]
can be considered as an element $\lan \und{\D}_X\ran$ in the Fock space associated with the adelic symplectic space.

We define the {\em adelic map} $\und{\ }: \K^{\infty} \to \und{\K}^{\infty}$ by
\[ \f = (\f_1,\f_2,\dots,\f_r,\dots) \mapsto \und{\f} = \{ f_r^{(\zeta)} \}: \
f_r^{(\zeta)}:=\Psi^r(\f_r(q^{1/m}/\zeta)) ,\]
where the last expression is to be expanded into a Laurent series near $q=1$ after applying Adams' operations $\Psi^r$, acting naturally on $K=K^0(X)\otimes \Lambda$, and by $\Psi^r(q)=q^r$ on functions of $q$. The residue theorem implies that the adelic map is symplectic:
\begin{align*} \Omega^{\infty} (\f,\g) &= \sum_{r=1}^{\infty}\frac{\Psi^r}{r}\sum_{\zeta} \Res_{q=\zeta} (\f_r(q^{-1}),\g_r(q))\, \frac{dq}{q}\\
  &= \sum_{r=1}^{\infty}\frac{1}{r^2}\sum_{\zeta}  \Res_{q=1}(\Psi^r(\f_r(q^{-1/m}\zeta)),\Psi^r(\g_r(q^{1/m}/\zeta))^{(r)} \frac{dq^{r/m}}{q^{r/m}} \\
  &= \sum_{\zeta}\frac{1}{m(\zeta)}
  \sum_{r=1}^{\infty}
  \frac{1}{r} \Res_{q=1}\ (\und{\f}_r^{(\zeta)}(q^{-1}), \und{\g}^{(\zeta^{-1})}_r(q))^{(r)} \frac{dq}{q} .\end{align*}
Our ``higher genus quantum RR formula'' can be stated this way.

\vspace{3mm}

{\tt Main Theorem.} {\em The adelic map $\und{\ }: (\K^{\infty}, \Omega^{\infty}) \to (\und{\K}^{\infty}, \und{\Omega}^{\infty})$ between the symplectic loop spaces transforms the adelic quantum state $\lan \und{D}_X\ran $ into the total descendant potential $\lan \D_X\ran $ of permutation-equivariant quantum K-theory of the target K\"ahler manifold $X$.}

\vspace{3mm}

How does a map between symplectic spaces map respective Fock spaces? Elements of the Fock space are functions on the symplectic space constant in the direction of the negative space of a chosen Lagrangian polarization. A map between symplectic spaces respecting the negative spaces of the chosen polarizations induces a map between the quotients, and hence maps the Fock spaces naturally (in the reverse direction). When the given polarizations disagree, one needs first to change one of them to identify the models of the Fock space based on different polarizations by the construction of Stone-von Neumann's theorem, and only after that apply the natural pull-back.

In the situation of our theorem, the polarizations disagree, and the precursory change of polarization in the adelic space is one of the key ingredients of the relation between $\D_X$ and $\und{\D}_X$ as generating functions.

The space $\K^{\infty}_{-}$ consists of sequences
$\f = (\f_1,\f_2,\dots,\f_r,\dots)$ of vector-values rational functions of $q$
with poles at roots of unity $\zeta$, but vanishing at $q=\infty$ and having
no pole at $q=0$. Such rational functions uniquely decompose into the sums of
their partial fractions, $\f_r = \sum_{\zeta} \f_r^{(\zeta)}$, i.e. reduced rational functions of $q$ with only one pole $q=\zeta$. In fact the negative space of polarization $\und{\K}^{\infty}_{-}$ in the adelic space (we've neglected to describe it so far, but it is involved in the interpretation of the infinite product $\und{\D}_X$ as an element of the Fock space) is exactly the direct sum of subspaces $\{ \Psi^r(\f_r^{(\zeta)}(q^{1/m(\zeta)}/\zeta)) \} \subset \K_r^{(\zeta)}=K((q-1))$ obtained from such partial fractions.

\pagebreak

By the way, we encounter here an interesting phenomenon impossible in finite-dimensional symplectic geometry. The adelic map $\und{\ }: \K^{\infty} \to \und{\K}^{\infty}$ is a symplectic injection which embeds the Lagrangian subspace $\K^{\infty}_{+}$ into the much bigger Lagrangian subspace $\und{\K}^{\infty}_{+}$, but it identifies the Lagrangian subspaces $\K_{-}^{\infty}$ and $\und{\K}_{-}^{\infty}$ considered as quotient spaces $\K^{\infty}/\K_{+}^{\infty}$ and $\und{\K}^{\infty}/\und{\K}_{+}^{\infty}$. 

At the same time, the image of $\K_{-}^{\infty}$ under the adelic map does not coincide with $\und{\K}_{-}^{\infty}$, and it is now easy to understand why:
the image of $\f_r^{(\zeta)}$ consists of the expansions of
$\Psi^r(\f_r^{(\zeta)}(q^{1/m(\eta)}/\eta))$ for all roots of unity $\eta$, and
not only for $\eta=\zeta$ where the partial fraction $\f_r^{(\zeta)}$ has its pole. Consequently, the relation between the quantum states $\lan \D_X\ran $ and $\lan \und{\D}_X\ran $ described in the theorem actually means that the total descendant potential $\D_X$ is obtained from the infinite product $\und{\D}_X$ as
\vspace{-2mm}
\[ \D_X = \text{pull-back by}\ \und{\ }: \K_{+}^{\infty}\subset \und{\K}_{+}^{\infty}\ \text{of}\ \ e^{ \frac{1}{2}\sum_r r\Psi^r(\h \sum \nabla_{\eta, \zeta})} \otimes_{M=1}^{\infty}\D_{X/\ZZ_M}^{tw}. \]
\vspace{-0.2mm}
Here $\nabla_{\eta,\zeta}$ are certain 2nd order differential operators
whose coefficients are tautologically determined by expansions of partial fractions with poles at roots of unity $\zeta$ into power series near all other
roots of unity, while the embedding $\und{\ }: \K_{+}^{\infty}\to \und{\K}_{+}^{\infty}$  maps sequences $\t=(\t_1,\t_2,\dots,\t_r\dots)$ of Laurent polynomial $\t_r \in K[q,q^{-1}]$ into the collection of power series expansions
$\Psi^r (\t_r(q^{1/m(\zeta)}/\zeta))$ of the Laurent polynomials at the roots of
unity.

The above description of our main formula is neither complete not totally accurate, and should be supplemented with further clarifications.

1. The quantum state $\lan \D_X \ran$ differs from the total descendant potential $\D_X$ (though both are functions on $\K^{\infty}/\K_{-}^{\infty}=\K_{+}^{\infty}$)
by the translation of the origin called the {\em dilaton shift}:
$\lan \D_X \ran (\vv+\t) = \D_X (\t)$, where $\vv = ((1-q)\1, (1-q)\1, \dots)$,
and $\1$ stands for the unit element in $K^0(X)$. Likewise,
$\lan \D_{X/\ZZ_M}^{tw} \ran ((1-q)\1+t) = \D_{X/\ZZ_M}^{tw}(t)$. Here $\1$ belongs to the unit sector, i.e. among the components $t^{(\zeta)}\in K[[q-1]]$ labeled by the $M$th roots of unity $\zeta$ only the component with $\zeta=1$ is dilaton-shifted.

2. In the generating functions for GW-invariants, one weighs contributions of degree-$d$ stable maps by the binomials $Q^d$ in Novikov's variables $Q=(Q_1,\dots,Q_r)$, where $r=\operatorname{rk} H_2(X,\ZZ)$. Novikov's variables are adjoined to the ground $\lambda$-ring $\Lambda$ so that $\Psi^r Q^d := Q^{rd}$. Furthermore, the expression ``rational functions'' (``Laurent series,'' ``power series'', etc.) of $q$ should be understood as
formal $Q$-series whose coefficients are rational functions (formal Laurent series, power series etc.) of $q$, and the notations like $K[q,q^{-1}]$, $K((q-1))$, etc. have to be understood in the sense of such a $Q$-adic completion.

3. To avoid some divergences, we require that $\Lambda$ is a local algebra with the maximal ideal $\Lambda_{+}$, that Adams' operations
respect the filtration by its powers: $\Psi^r\Lambda_{+}\subset \Lambda_{+}^r$, and assume that the components of the variables in generating functions lie in $\Lambda_{+}$. In particular, the quantum states $\lan \D_X\ran$, $\lan \und{\D}_X\ran$, etc. are functions on $\K_{+}^{\infty}$, $\und{\K}_{+}^{\infty}$, etc. defined in a $\Lambda_{+}$-neighborhood of the dilaton shift. 

4. A peculiar phenomenon overlooked in the previous discussion is
that the symplectic structure $\Omega^{\infty}$, the adelic map, and other  ingredient of our formalism are not $\Lambda$-linear in the usual sense. For instance, for $\nu \in \Lambda$ and $\f = (\f_1,\f_2,\dots,\f_r,\dots ) \in \K^{\infty}$,
the adelic image $\und{\nu \f} = (\nu \und{\f}_1, \Psi^2(\nu) \und{\f}_2, \dots, \Psi^r(\nu) \und{\f}_r,\dots )$,
i.e. the map between the $r$th components is {\em linear relative to} the scalar transformation $\Psi^r$.\footnote{Perhaps one can rectify this by noticing that {\em de facto} $\D_X$ depends not on $\t_r$, $r=1,2,3,\dots$, but on $\Psi^r\t_r$.}     

5. The previous feature manifests in the quantization formalism as well. Namely
the Planck constant, which needs to be adjoint to the ground ring $\Lambda$, is
acted upon by Adams' operations as $\Psi^r\h := \h^r$.  Respectively, $\h^r$ plays the role of the Planck constant in the quantization formalism on the $r$th component of the adelic space $\und{\K}^{\infty}$. This is manifest in our
formula $\sum_r r\h^r \Psi^r\sum \nabla_{\eta,\zeta}$ for the propagator, where
$\nabla_{\eta,\zeta}$ are 2nd order differential operators. 

6. This brings up the question about the status of the Planck constant in the
adelic product $\otimes_{M=1}^{\infty} \D_{X/\ZZ_M}^{tw}$ since each factor mixes
up sectors with different values of the index $r$. In fact the quantum state
$\lan \D_{X/\ZZ_M}^{tw} \ran (t, \h, Q))$ (i.e. the generating function
for twisted fake K-theoretic GW-invariants of the orbifold $E_M$ {\em after} the dilaton shift) is homogeneous (due to the so-called dilaton equation):
\[  \lan \D_{X/\ZZ_M}^{tw} \ran (t, \h, Q)) = \h^{\frac{M \dim K^0(X)}{48}} \lan \D_{X/\ZZ_M}^{tw}\ran (\frac{t}{\sqrt{\h}}, 1, Q) .\]
By the rules of quantum mechanics, scalar factors don't affect ``quantum states.'' The accurate definition of the infinite tensor product in our main theorem is
\[ \lan \und{\D}_X \ran \, \, (\{ t_r^{(\zeta)}\}, \h, Q) = \otimes_{M=1}^{\infty} \lan \D_{X/\ZZ_M}^{tw}\ran \left(\frac{\{t_{M/m(\zeta)}^{(\zeta)}\} }{\sqrt{\h}}, 1, Q^M\right) .\]
Note the change of $Q$ into $Q^M$ in the $M$th factor.

7. Our final remark here is about equivariant generalizations of the theorem. In applications of GW-theory, the target space is often equipped with an action
of a torus $T$, and all holomorphic Euler characteristics are replaced with the characters of the $T$-action on the sheaf cohomology. In particular,  Lefschetz' fixed point localization technique, when combined with the formalism of symplectic loop spaces, leads to dealing with fractions
of the form $1/(1-q^m\tau)$, where $\tau$ is a coordinate on $T$, and the poles in $q$ are at roots of $1/\tau$ rather than roots of unity. Nevertheless our theory carries over {\em verbatim} to the equivariant case. Namely, the homotopy theory construction of equivariant K-theory yields $K^0_T(pt)=K^0(BT)$ which is not the character ring of $T$, but its completion into functions on $T$ defined in the formal neighborhood of the identity. Our ground $\lambda$-algebra $\Lambda$ should be changed into $\Lambda \otimes K^0_T(pt)$. To make sense, the above fractions must be expanded into series in $\tau-1$ with coefficients in rational fractions of $q$
having poles at roots of unity only:
\[ \frac{1}{1-q^m\tau}=\frac{1}{1-q^m-q^m(\tau-1)}=
\sum_{n=0}^{\infty}\frac{q^{mn}(\tau-1)^n}{(1-q^m)^{n+1}}.\]
Thus, in the homotopy theory interpretation of $T$-equivariant K-theory,  localization to fixed points of $T$ makes no sense, but 
our ``quantum RR formula'' holds unchanged for $T$-equivariant GW-invariants, which take values in $\Lambda\otimes K^0_T(pt)$.

\section{Redefining the invariants}

Let us recall and generalize the definition of permutation-equivariant K-theoretic GW-invariants given in Part I, and of the mixed genus-$g$ potential given in Part VII.

Let $X$ be a compact K\"ahler manifold, $K:=K^0(X)\otimes \Lambda$, where  $\Lambda$ is a local $\lambda$-algebra that contains Novikov's ring as it was  explained in Introduction. 

Let $X_{g,n,d}$ be the moduli space of degree-$d$ stable maps to $X$ of complex curved of arithmetic genus $g$ with $n$ marked points, and let $h\in S_n$ be
a permutation, acting on the moduli space by renumbering the marked points.
Let $V$ be a holomorphic vector bundle over $X_{g,n,d}$ equivariant with respect to the action of the permutation $h$. Then the sheaf cohomology
$\pi_*(V):=H^*(X_{g,n,d}; V\otimes \O_{g,n,d})$, where $\O_{g,n,d}$ is the ($S_n$-invariant) virtual structure sheaf introduced by Y.-P. Lee \cite{YPLee}, inherits the action of $h$. Therefore the supertrace $\str_h \pi_*(V)$ is defined.

Denote $l_k=l_k(h)$ the number of cycles of length $k$ in the cycle structure of $h$, and by $\ll=(l_1,l_2,l_3,\dots)$ the corresponding partition of $n=\sum rl_r$.  Our current goal is to define {\em correlators} 
\[  \lan \uu_1,\dots, \uu_{l_1}; \vv_1, \dots, \vv_{l_2}; \dots ; \ww_1, \dots, \ww_{l_r}; \dots \ran_{g,\ll,d},\]
where the {\em inputs} $\uu_i, \vv_j, \ww_k, \dots $ are elements of $K^0(X)\otimes \Lambda [q,q^{-1}]$. Note that groups of the seats in the correlator have lengths $l_1$, $l_2$ etc., and the total number $\sum l_r$ of the seats is equal to the number of {\em non-empty} cycles.

Let $\gs_1,\dots,\gs_r$ be indices of the marked points cyclically permuted by $h$, and let out of all the $l_r$ cycles of length $r$, this be the $k$th cycle. We take the $h$-equivariant bundle $W_k$ on $X_{g,n,d}$ determined by the input $\ww_k = \sum_m \phi_m q^m$ ($\phi_m \in K^0(X)$) in the form
\[ W_k := \bigotimes_{\a=1}^r \sum_m (\ev_{\gs_\a}^*\phi_m) L_{\gs_\a}^m,\]
where $\ev_{\gs_{\a}}: X_{g,n,d}\to X$ is the evaluation map, and $L_{\gs_\a}$ is the universal cotangent line bundle at the marked point with the index $\gs_\a$.
This way, for each cycle of length $1$, $2$, etc. we associate the inputs $\uu_i$, $\vv_j$, etc. and define respectively the bundles $U_i$, $V_j$, etc.
We define the above correlator as  
\[ \prod_{r=1,2,\dots} r^{-l_r}\str_h H^*\left(X_{g,n,d}; \O_{g,n,d}\bigotimes_{i=1}^{l_1} U_i \bigotimes_{j=1}^{l_2} V_j \cdots \bigotimes_{k=1}^{l_r}W_k \cdots \right) .\]   
The factor in front of the supertrace is motivated by the number
$n!/\prod_r r^{l_r} l_r!$ of permutations with the cycle structure described by the partition $\ll$.

Note that the correlator is poly-additive with respect to each input. Namely,
if $\ww_k = \ww_k'+\ww_k''$, then
\[ \bigotimes_{\a=1}^r \ww_k(L_{\a}) = \sum_{I\subset \{\ 1,\dots,n\} } \bigotimes_{\a \in I} \ww'_k(L_{\a})\bigotimes_{\b\notin I} \ww''(L_{\b}).\]
The sheaf cohomology splits into $2^r$ summands accordingly, but the summands with $I\neq \emptyset$ or $\{ 1,\dots, n\}$ are permuted by $h$ non-trivially, and hence don't contribute to $\str_h$. Therefore
\[ \lan \dots, \ww_k,\dots\ran_{g,\ll,d}= \lan \dots, \ww'_k,\dots\ran_{g,\ll,d}+\lan \dots, \ww''_k,\dots\ran_{g,\ll,d}.\]

We extend the correlator to inputs from $\K_{+}:=K^0(X)\otimes \Lambda [q,q^{-1}]$ in the way {\em linear relative to $\Psi^r$} on each input corresponding to the cycles of length $r$, i.e.    
\[  \lan \dots, \nu \ww_k,\dots\ran_{g,\ll,d} = \Psi^r(\nu) \lan \dots, \ww_k,\dots\ran_{g,\ll,d}.\]  
This is motivated by the fact that if $\Lambda = K^0(Y)$, then for a vector bundle $\nu$ on $Y$, the trace bundle of the cyclic permutation of the factors in $\nu^{\otimes r}$ coincides with $\Psi^r(\nu)$.

Now, we define the genus-$g$ potential of permutation-equivariant quantum K-theory of $X$ as the sum over degrees and partitions $\ll$ of all $n=0,1,2,\dots$:
\[ \F_g(\t) = \sum_d Q^d \sum_{\ll} \frac{1}{\prod_r l_r!} \lan \dots \t_1 \dots ; \dots, \t_2, \dots ; \dots \ran_{g,\ll,d}.\]
Here $\tt=(\t_1,\t_2,\dots, \t_r,\dots)$, each $\t_r\in \K_{+}$,
and all the inputs in the correlator corresponding to the cycles of length $r$
are taken to be the same and equal $\t_r$.

\medskip

{\tt Remark.} The correlators $\lan \uu, \dots, \uu \ran_{g,n,d}^{S_n}$ defined in Part I by taking averages over $S_n$ can be expressed in terms of the above correlators via re-summation over the conjugacy classes labeled by partitions $\ll$ of $n$:
\[ \lan \uu, \dots, \uu \ran_{g,n,d}^{S_n}=\frac{1}{n!}\sum_{h\in S_n}
  \str_h [\uu,\dots,\uu]_{g,n,d} =\sum_{\ll} \frac{1}{\prod_r l_r!} \lan \uu ; \dots ; \uu \ran_{g,\ll,d}.\]
  Respectively the mixed genus-$g$ potential of Part VII 
  \[ \sum_{m,n\geq 0, d} Q^d \lan \x,\dots, \x; \t,\dots,\t\ran_{g,m+n,d}^{S_n} = \sum_d Q^d \sum_{\ll} \frac{1}{\prod_r l_r!} \lan \x + \t; \t ; \t; \dots \ran_{g,\ll,d}\]
coincides with the specialization of $\F_g$ to the inputs $\t_1=\x+\t$,
$\t_2= \t, \t_3=\t, \dots$. 

\medskip

While moduli spaces $X_{g,n,d}$ parameterize stable maps of {\em connected} curves, the {\em total} descendant potential is to account for contributions of possibly disconnected curves, as well as for symmetries of such curves caused by permutations of identical connected components. 

Abstractly speaking, if $\nu \in \gL$ represents the contribution of ``connected'' objects, then the sum over $n$ of contributions of objects with $n$ components is given by
\[ \sum_{n\geq 0} \frac{1}{n!} \sum_{h\in S_n} \prod_{k>0} \Psi^k(\nu)^{l_k(h)} = \sum_{\ll}\prod_{k>0}\frac{(\Psi^k(\nu)/k)^{l_k}}{l_k!} = e^{\sum_{k>0} \Psi^k(\nu)/k}. \]
This motivates the following definition of
the {\em total descendant potential} of the permutation-equivariant quantum K-theory on $X$:
\[ \D_X := e^{\textstyle \sum_{g\geq 0} \left[\sum_{k>0}\h^{k(g-1)} \Psi^k (R_k \F_g)/ k\right] },\] 
where
$(R_k \F) (\t_1,\t_2,\dots,\t_r,\dots) :=\F(\t_k,\t_{2k},\dots,\t_{rk},\dots)$. 

In order to explain the rescaling $R_k$ of the indices in the variables $\t_r$, note that automorphisms of $\bigsqcup \prod_{\a=1}^k X_{g_\a, n_\a, d_\a}$ induced by cyclic permutations
of $k$ connected components of a disconnected curve accompanied by a renumbering $h$ of marked points, generate traceless operators on the sheaf cohomology unless $g_\a, n_\a, d_\a$ don't depend on $\a$, and $h^k$ renumbers
the marked points of all components separately in consistent ways. In this case, we have an automorphism of $X_{g,n,d}^k$ whose $k$th power is the automorphism of each factor $X_{g,n,d}$ induced by the renumbering $h^k$. If the orbit of one of the marked points under the renumbering $h^k$ has order $r$, then the orbit under the renumbering $h$ has order $rk$. Therefore the input corresponding to this cycle of marked points must be $\t_{rk}$.

Finally, the factor $\h^{k(g-1)}$, whose exponent is $-1/2$ times the Euler
characteristic of $k$ copies of a genus-$g$ Riemann surface, can be interpreted
as $\Psi^k(\h^{g-1})$ by adjoining $\h$ to $\Lambda$ and setting $\Psi^k(\h)=\h^k$. 

\vspace{3mm}

Note that all $\F_g$ can be recovered from ${\mathcal G}:=\log \D_X$ by M\"obius' exclusion-inclusion formula
\[ \sum_g \h^{g-1}\F_g = \prod_{p\ \text{prime}}\left(1-\frac{\Psi^p}{p}\, R_p \right)\ {\mathcal G}.\]

\section{Kawasaki's Riemann--Roch formula}

The expression of K-theoretic GW-invariants in terms of cohomological ones is based on the use of the virtual variant \cite{ToK} of Kawasaki's Riemann--Roch formula \cite{Kaw}.

Let $\M$ be a compact complex orbifold, and $V$ be a holomorphic orbibundle on $\M$. The holomorphic Euler characteristic of $V$, defined in terms of \v{C}ech cohomology as $\chi (\M; V) :=\sum_i (-1)^i \dim H^i(\M; V)$, is expressed by Kawasaki's RR formula in cohomological terms of the {\em inertia orbifold} $I\M$:
\[ \chi (\M; V) = \chi^{fake} \left(I\M; \frac{\tr_h V}{\str_h \wedge^{\bullet} N^*_{I\M}}\right).\]
Recall that a point in $I\M$ is represented by a pair $(x, h)$ where $x\in \M$,
and $h \in \Gamma (x)$ is an element of the {\em inertia group} of $x\in \M$ (i.e. the group of local symmetries of $x$ in the orbifold structure).
In the formula, $N^*_{I\M}$ denotes the conormal bundle to the stratum of fixed points of the symmetry $h$. The bundle $V$ can be restricted to the stratum and decomposed into eigenbundles $V_{\lambda}$ of $h$ corresponding to the eigenvalues $\lambda$. The trace operation $\tr_h V$ denotes the virtual bundle $\sum_{\lambda} \lambda V_{\lambda}$, and the supertrace $\str_h$ in the denominator denotes the similar operation on the $\ZZ_2$-graded bundle $\wedge^{\bullet} N^*_{I\M}$. Finally, the notation $\chi^{fake}$ stands for the {\em fake} holomorphic Euler characteristic of an orbibundle over an orbifold:
\[ \chi^{fake} (M; W) := \int_M \ch (W) \td (T_M) ,\]
where $\ch (W)$ is the Chern character of the orbibundle $W$, and $\td (T_M)$ is the Todd class of tangent orbibundle $T_M$ (both defined over $\QQ$). 

In effect, the RHS of Kawasaki's RR formula is the sum of certain fake holomorphic Euler characteristics, i.e. of certain integrals over the strata
of the inertia orbifold, which are rational numbers adding up to the integer defined by the LHS.

\medskip

It is no accident that Kawasaki's RR formula resembles Lefschetz' holomorphic fixed point formula. To make the connection, let $h$ be an automorphism of a holomorphic bundle $V$ over a compact complex manifold $\tilde{\M}$. For our goals it suffices to assume that $h$ belongs to a finite group $G$ of such automorphisms (although abstractly speaking this restriction can be relaxed). Lefschetz' fixed point formula computes the supertrace of $h$ on the sheaf cohomology as an integral over the fixed point submanifold $\tilde{\M}^h$:
\[ \str_h H^*(\tilde{\M}; V) = \chi^{fake} \left( \tilde{\M}^h ; \frac{\tr_h V}
{\str_h \wedge^{\bullet}N^*_{\tilde{\M}^h}}\right) .\]
On the other hand, $V$ can be considered as an orbibundle over the quotient orbifold $\M:=\tilde{\M}/G$, and the holomorphic Euler characteristic $\chi (\M; V)$ of the orbibundle can be found as the average over $G$:
\[ \frac{1}{|G|}\sum_{h\in G} \str_h H^*(\tilde{\M}; V) =
\frac{1}{|G|}\sum_{h\in G}\chi^{fake} \left( \tilde{\M}^h ; \frac{\tr_h V}
{\str_h \wedge^{\bullet}N^*_{\tilde{\M}^h}}\right) .\]
The last sum coincides with the right hand side of Kawasaki's RR formula
on $\M=\tilde{\M}/G$ since in the global quotient case
\[ I\M = \left[ \bigsqcup_{h\in G} \tilde{\M}^h \right] / G .\] 

In fact, we need a combination of Kawasaki's RR with Lefschetz' fixed point formula, computing $\str_h H^*(\M ; V)$ where $h$ is a finite order automorphism of an {\em orbi}bundle $V$ over an {\em orbi}fold $\M$: 
\[ \str_h H^*(\M; V) = \chi^{fake} \left( I\M^h; \frac{\tr_{\tilde{h}} V}{\str_{\tilde{h}} \wedge^{\bullet}N^*_{I\M^h}} \right) ,\]
where the ``fixed point inertia orbifold'' $I\M^h$ can be described as follows.
Let $x\in \M$ be a fixed point of $h$, and $U_x \to U_x/\Gamma(x) \subset \M$
be its orbifold chart. The transformation $h$ can be lifted to automorphisms $\tilde{h}$ of the chart (and of the bundle $V$ over the chart) in $|\Gamma(x)|$ possible ways. Each transformation $\tilde{h}$ has a fixed point
submanifold $U_x^{\tilde{h}} \subset U_x$ whose union is $\Gamma(x)$ invariant.
The quotient $\left[ \bigcup_{\tilde{h}} U_x^{\tilde{h}} \right]/ \Gamma(x)$ provides the local description of the orbifold $I\M^h$ near $x \in \M$. The ingredients $\tr_{\tilde{h}} V$ and $N_{I\M^h}$ are obtained from the fibers of $V$ over $U_x^{\tilde{h}}$ and from the normal space to $U_x^{\tilde{h}}$ in $U_x$
respectively.

A justification of Lefschetz-Kawasaki's RR formula can be obtained formally from Kawasaki's RR formula applied to the orbifold $\M/G$ where $G$ is the cyclic group generated by $h$. Indeed, let $\CC_{\lambda}$ denotes the $1$-dimensional representation of $G$ where $h$ acts by a root of unity $\lambda$. Then
\[ \str_h H^*(\M; V) := \sum_{\lambda} \lambda\ H^*(\M; V)_{\lambda} =
\sum_{\lambda} \lambda \ \chi (\M/G; V \otimes \CC_{\lambda^{-1}}) .\]
The last sum can be computed on the inertia orbifold $I(\M/G)$ using Kawasaki's RR. However $\sum_{\lambda} \lambda \, \CC_{\lambda^{-1}}$ is a virtual representation of $G$ whose character equals $|G|$ on $h$ and equals $0$ on all other elements of $G$. Therefore only the strata of $I(\M/G)$ made of fixed points of $\tilde{h}$ will contribute. Note that the factor $|G|$ from the character is compensated by the factor $1/|G|$ arising from the comparison between the fundamental classes of strata in $I\M^h$ with those in $I(\M/G)$.

\smallskip

In applications to quantum K-theory, the orbifold $\M$ is replaced with moduli space $X_{g,n,d}$ of stable maps to $X$, which are virtual orbifolds, or with
products of such spaces (since the curves are allowed to be disconnected).
An automorphism $h$ of such a product is induced by a renumbering of the marked points on the curve. A fixed point of $h$ is represented by a stable map $\phi: \Sigma \to X$ for which there exists a symmetry accomplishing the required permutation $h$, i.e. there exists an isomorphism $\tilde{h}: \Sigma \to \Sigma$ which permutes the marked points by $h$, and such that $\phi \tilde{h} = \phi$. It is the result of \cite{ToK} which justifies the application of Kawasaki's RR to {\em virtual} orbifolds.\footnote{The set-up of the virtual Kawasaki RR is axiomatic, but it eventually employs Kawasaki's RR theorem for (ambient) compact orbifolds. For moduli spaces of stable maps,
the existence of such ambient orbifolds is easily obtained in genus 0 by
projective embedding of $X$ (since $\M_{0,n}(\CC P^n, d)$ are orbifolds). In higher genus, the existence of such {\em compact} ambient orbifolds is a result of A. Kresch \cite{Kr}. Of course, one expects Kawasaki's RR formula to remain true for compactly supported orbisheaves on non-compact orbifolds (which would settle this technical issue in a more natural way). For compactly suppotred sheaves on manifolds, this was proved in \cite{OTT} some quarter of a century later than Hirzebruch's celebrated result for compact manifolds. The orbifold story develops slower, and almost 40 years after Kawasaki's result \cite{Kaw}, its vision for compactly supported orbisheaves seems still missing in the literature. The most promising approximations we could find were \cite{Far} and \cite{Ed}.}     
Respectively, our generating function $\D_X$ (which incorporates contributions of all stable maps and all renumberings of the markings) can be described in terms of suitable  fake holomorphic Euler characteristics on the strata of the inertia orbifold $I\M^h$. We will call them {\em Kawasaki strata}. They parameterize {\em stable maps with prescribed symmetries}, i.e. equivalence classes of pairs $(\phi, \tilde{h})$, where $\phi$ is a stable map of a (possibly disconnected) curve to $X$, and $\tilde{h}$ is a symmetry of the map, accomplishing a (possibly non-trivial) permutation of the marked points.

\begin{figure}[htb]
\begin{center}
\epsfig{file=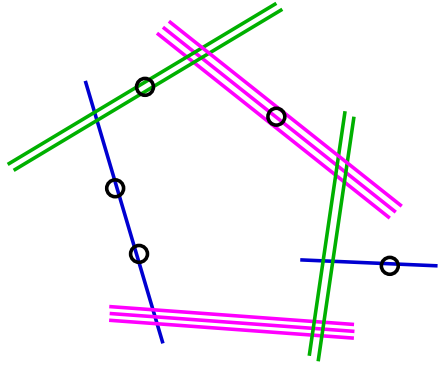} 
\end{center}
\caption{Stable maps with prescribed symmetries}
\end{figure}

How does a Kawasaki stratum look like? Given a stable map $\phi: \Sigma \to X$
with a symmetry $h$ (note that now on we omit the tilde), it defines the map of the quotient $\hat{\Sigma}$ of the curve $\Sigma$ by the cyclic group generated by $h$. On Figure 1, we attempt to
show a typical picture of a (connected) quotient curve. The quotient map
$\Sigma \to \hat{\Sigma}$ may have different number of branches (shown as the multiplicity of lines) over different irreducible components of $\hat{\Sigma}$. This shows that the summation over Kawasaki strata will have the structure of Wick's formula of {\em summation over graphs}. The vertices of the graphs represent contribution of Kawasaki strata parameterizing irreducible quotient maps, while the edges correspond to the nodes connecting the irreducible components.

Furthermore, an $M$-fold quotient map $\Sigma \to \hat{\Sigma}$ over an irreducible curve $\hat{\Sigma}$ can be described as the principal $\ZZ_M$-bundle over the complement to marked and nodal points, possibly ramified at such points. Consequently, Kawasaki strata representing the
vertices can be identified with {\em moduli spaces of stable maps to the orbifold target spaces $X/\ZZ_M$} ($=X \times B\ZZ_M$ in the notation of \cite{JK}, i.e. assuming that $\ZZ_M$ acts trivially on $X$).

We will denote by $\D_{X/\ZZ_M}^{fake}$the total descendant potential of the {\em fake} quantum K-theory of the orbifold $X/\ZZ_M$. Using the results \cite{ToTs},
one can obtain the K-theoretic counterpart to the theorem of Jarvis-Kimura \cite{JK} and express $\D_{X/\ZZ_M}^{fake}$ in terms of $\D_X^{fake}$, the total descendant potential of quantum K-theory of $X$. The latter can be, in its turn, expressed in terms of the cohomological total GW-potential $\D_X^H$, using the quantum Hirzebruch-Riemann-Roch formula \cite{Co, CGL} for fake GW-invariants with values in complex cobordisms, specialized to the case of complex K-theory. However, the vertex contributions in our Wick's formula are not $\D_{X/\ZZ_M}^{fake}$, but some {\em twisted} fake K-theoretic GW-invariants of
these orbifolds. This means that the virtual fundamental classes of moduli spaces of stable maps to $X/\ZZ_M$ need to be systematically modified --- in fact by the factors accounting for the denominators in the Kawasaki-RR formula. The total descendant potential $\D_{X/\ZZ_M}^{tw}$ for suitably twisted
fake quantum K-theory of $X/\ZZ_M$ can be expressed in terms of $\D_{X/\ZZ_M}^{fake}$ using the results of Tseng \cite{Ts} and Tonita \cite{ToT}.

In the next two sections, we first explain (or recall) how to pass from $\D_X^H$ to $\D_X^{fake}$, and then to $\D^{fake}_{X/\ZZ_M}$. Then we will formulate the twisting result relating $\D_{X/\ZZ_M}^{fake}$ with $\D_{X/\ZZ_M}^{tw}$. Then the vertex contributions of our graph summation formula will be described, roughly speaking, as the product $\bigotimes_{M=1}^{\infty} \D_{X/\ZZ_M}^{tw}$ over all $M=1,2,3, \dots$, leading to the concise quantum-mechanical description of $\D_X$ given in Introduction.          

\section{Symplectic loop spaces and quantization}

The formalism of symplectic loop spaces and their quantizations starts with the datum: a vector {\em space} $H$ (or a module over a ground ring $\Lambda$), a symmetric $\Lambda$-valued {\em Poincar{\'e} pairing} $(\cdot, \cdot)$ on $H$, and a nonzero vector $v\in H$. Using this datum, one cooks up a {\em loop space} $\H$, equipped with a symplectic $\Lambda$-valued form $\Omega$, a Lagrangian {\em polarization} $\H:=\H_{+}\oplus \H_{-}$, and a vector $\vv \in \H_{+}$ called the {\em dilaton shift}.

Given a sequence of functions $\F_g: \H_{+}\to \Lambda$, one combines them into the {\em total descendant potential} $\D:=e^{\textstyle \sum \h^{g-1}\F_g}$, and interprets the latter as an ``asymptotical element'' in the Fock space associated with $(\H, \Omega)$ by  
lifting from $\H_{+}$ to $\H$ the dilaton-shifted function $\t \mapsto \D(\t-\vv)$ from $\H_{+}$ so that it stays constant along the Lagrangian subspaces parallel to $\H_{-}$.    

According to the ideology of quantum mechanics, the Heisenberg Lie algebra of
the symplectic space acts irreducibly in the Fock space (of functions constant in the direction of $\H_{-}$), which by Schur's lemma, projectively identifies Fock spaces defined using different polarizations. Furthermore, the symplectic group moves the polarizations around, which therefore defines a projective action of the Lie algebra of quadratic hamiltonians on the Fock space. Explicit formulas for this action provide the standard
quantization of quadratic hamiltonians. Namely, let $\{\ q_\a \}$ be coordinates on $\H_{+}$, and $\{ p_\a \}$ the Darboux-dual coordinates on
$\H_{-}$. Then the quantization $\widehat{\ }$ of Darboux monomials is given by the multiplication and differentiation operators on functions of $\{ q_\a\}$:
  \[ \widehat{q_\a q_\b} := \h^{-1} q_\a q_\b, \  \ \widehat{q_\a p_\b} :=
  q_{\a} \partial_{q_{\b}},\ \ \widehat{p_\a p_\b}:= \h\,  \partial_{q_\a}\partial_{q_\b}.\] 
  Finally, given a linear symplectic transformation $\square$ on $(\H, \Omega)$, the {\em Stone - von Neumann quantization} of it acts on the Fock space by the operator $\widehat{\square} := e^{\widehat{\log \square}}$.

  A typical application of this formalism in GW-theory relates generating functions for two kinds of GW-invariants as follows. The functions $\D^i$, $i=1,2$, are lifted to asymptotical elements $\lan \D^i\ran$ of
the respective Fock spaces associated with symplectic loop spaces $(\H^i, \Omega^i)$ using Lagrangian polarizations $\H_{\pm}^i$ and dilaton shifts $\vv_i$. The respective quantum states are related by
  \[ \lan \D^1 \ran = \widehat{\qch}\ \widehat{\square}\ \lan \D^2 \ran ,\]
  where $\square$ is a suitable symplectic automorphism of $(\H^2,\Omega^2)$,
  while the ``quantum Chern character'' $\qch: \H^1\to \H^2$ is a symplectic isomorphism (i.e. $\qch^*\Omega^2=\Omega^1$), and hence identifies the respective Fock spaces. Note that the isomorphism $\qch$ may not respect the polarizations (in practice, $\qch$ respects $\H^i_{+}$, but not $\H^i_{-}$),
nor the dilaton shifts ($\qch \vv_1\neq \vv_2$). Consequently, the generating functions $\D^1$ and $\D^2$ are obtained from each other by three consecutive transformations: the quantized operator $\square$, the change of polarization, and the correction for the discrepancy in the dilaton shifts.

\medskip

To begin with cohomological GW-invariants of $X$, we set
\[ H:=H^{even}(X; \Lambda),\ \ (a,b)^H:=\int_X ab,\ \ v=\1,\]
take $\H$ to be the space $H((z))$ of Laurent series in one indeterminate $z$ with vector coefficients from $H$. We assume that the ground ring $\Lambda$ contains Novikov's variables, $Q$, and the Laurent series are {\em $Q$-adically convergent} for $z\neq 0$, i.e. that modulo any fixed power of $(Q)$, the series in question contain finitely many negative powers of $z$. We equip $\H$ with the symplectic form
\[ \Omega^H(\f,\g):= \Res_{z=0} (\f(-z), \g(z))^H\, dz,\]
and Lagrangian polarization $\H=\H_{+}\oplus \H_{-}$, where $\H_{+}$ consists of the power series part of the Laurent series, and $\H_{-}$ of their principal parts.

Recall that genus-$g$ generating functions for GW-invariants of $X$ are defined by
\[ \F_g^H(\t):=\sum_{d,n} \frac{Q^d}{n!} \int_{[X_{g,n,d}]} \prod_{i=1}^n \sum_{k=0}^{\infty} \sum_{\a} t_{k,\a} \ev_i^*(\phi_{\a}) \psi_i^k,\]
where $[X_{g,n,d}]$ is the virtual fundamental classes of the moduli spaces of stable maps to $X$, $\psi_i := c_1(L_i)$ is the 1st Chern class of universal cotangent line bundle at the $i$th marked point, and $\{ \phi_\a \}$ is a basis
in $H^{even}(X,\Lambda)$. They are functions of $\t = \sum_{k,\a} t_{k,\a} \phi_{\a} z^k$, which lie in $\H_{+}$. Respectively, the total descendant potential of
the cohomological GW-theory of $X$ is defined as $\D^H_X=e^{\sum_g \h^{g-1}\F_g^H(\t)}$, subject to the dilaton shift $\vv = -z\1$, i.e.
$\lan \D^H_X\ran (\t-z\1) = \D^H_X (\t)$. 

\medskip

In the fake quantum K-theory of $X$, one puts
\[ H:=K=K^0(X)\otimes \Lambda,\ (a,b):=\chi (X; a\otimes b) = \int_X \ch(a)\ch(b) \td(T_X),\]
uses $\K^{fake} = K((q-1))$, i.e. the space of $Q$-adically convergent Laurent series in $q-1$ with vector coefficients in $K$, and equips it with the symplectic form
\[ \Omega^{fake}(\f,\g):= \Res_{q=1} (\f(q^{-1}),\g(q))\, \frac{dq}{q},\]
and Lagrangian polarization $\K=\K_{+}\oplus \K_{-}$, taking $\K_{+}$ to consist of power series, and $\K_{-}$ of the principal parts of Laurent series in
$q-1$.

The genus-$g$ generating functions $\F_g^{fake}$ are defined on $\K_{+}$ by
\[ \F_g^{fake}(\t) = \sum_{d,n} \frac{Q^d}{n!} \chi^{fake}\left( X_{g,n,d};
\bigotimes_{i=1}^n \sum_{k=0}^{\infty}\sum_{\a} t_{k,\a} \ev_i^*(\phi_{\a}) (L_i-1)^k \right),\]
where $\{ \phi_\a \}$ form a basis in $K^0(X)$, and the fake holomorphic Euler characteristic of a bundle $V$ on $X_{g,n,d}$ is defined using the virtual fundamental cycle $[X_{g,n,d}]$ and the virtual tangent bundle bundle $T_{X_{g,n,d}}$:
\[ \chi^{fake}(X_{g,n,d}; V) := \int_{[X_{g,n,d}]} \ch(V) \td(T_{X_{g,n,d}}).\]
The total descendant potential of fake quantum K-theory is defined by
$\D_X^{fake}=e^{\sum \h^{g-1}\F_g^{fake}(\t)}$ as a function on $\K_{+}$ subject to the
dilaton shift by $\vv = (1-q)\1$, i.e. $\lan \D_X^{fake} \ran ((1-q)\1+\t)=\D_X^{fake}(\t)$. It is expressed in terms of $\D_X^H$ following \cite{Co, CGL}.

Namely, introduce the {\em quantum Chern character} $\qch: \K \to \H$ by 
\[ \K\ni \f =\sum_k f_k (q-1)^k \mapsto \sqrt{\td (T_x)} \sum_k \ch (f_k) (e^z-1)^k \in \H . \]
It is symplectic: $\qch^*\Omega^H = \Omega^{fake}$. Then
\[ \lan \D_X^{fake} \ran = \widehat{\qch}^* \ \widehat{\triangle}\
\lan \D_X^H \ran,\]
where $\triangle$ is the Euler--Maclaurin asymptotics of the infinite product
$\prod_{r=1}^{\infty} \td ((T_X-1)\otimes q^{-r})$. The equality holds up to a scalar factor explicitly described in \cite{Co}. Recall that the Euler--Maclaurin asymptotics of the product $\sqrt{S(E)}\prod_{r=1}^{\infty} S(E\otimes q^{-r})$, where $E$ is a vector bundle over $X$, $q$ is the universal line bundle (so that $c_1(q)=z$), and $S(\cdot)=e^{\textstyle\sum_k s_k\ch_k(\cdot)/k!}$ is an invertible multiplicaive characteristic class, is 
\[ e^{\textstyle{\sum_{m\geq 0} \sum_{l\geq 0} s_{2m-1+l}\frac{B_{2m}}{(2m)!}\ch_l(E)z^{2m-1}}},\]
  where $B_{2m}$ are Bernoulli numbers, and $\ch_l(E)$ in the exponent are understood as operators of classical multiplication in the cohomology algebra of $X$ by the components of the Chern character. 

  \medskip

Our next step is to describe in terms of $\D_X^{fake}$ the total descendant potential $\D^{fake}_{X/\ZZ_M}$ of the fake quantum K-theory of the orbifold $X/\ZZ_M$. The Grothendieck group $K^0(X/\ZZ_M)$ of orbibundles on $X/\ZZ_M$ is identified  with $K^0(X)\otimes \operatorname{Repr} (\ZZ_M)$. Respectively, the total descendant potential $\D_{X/\ZZ_M}^{fake}$ in the {\em fake} quantum K-theory of $X/\ZZ_M$ is a function on the space of vector power series 
\[ \t:=\sum_{\chi \in \operatorname{Repr}(\ZZ_M)} \t_{\chi} \chi, \]
where each $\t_{\chi}$ is a power series in $q-1$ with coefficients in $K^0(X) \otimes \gL$. In down-to-earth terms we have:
\[ \D^{fake}_{X/\ZZ_M}(\t)=\prod_{\chi \in\operatorname{Repr}(\ZZ_M) } \D^{fake}_X(\t_{\chi}).\]
This follows from the analogous cohomological result of Jarvis-Kimura \cite{JK}
by application of twisting theorems of Tseng \cite{Ts} and Tonita \cite{ToT}
(combined with a description of the virtual tangent bundles to the moduli spaces of stable maps to $X/\ZZ_M$). Alternatively, this result can be extracted from section 3 of their joint paper \cite{ToTs}.

To go on, we need to describe the element of the Fock space defined by $\D^{fake}_{X/\ZZ_M}$, and the respective symplectic loop space. We have
\[ H := K \otimes \operatorname{Repr}(\ZZ_M), \ \ (a,b):=\sum_{\chi} (a_\chi, b_\chi)^{fake}, \ \ v = \sum_{\chi} 1\chi.\]
Respectively the loop space 
\[ \K^{fake}_{X/\ZZ_M}=\K^{fake}_X\otimes \operatorname{Repr}(\ZZ_M),\]
is equipped with the symplectic form
\[ \Omega^{fake}_{X/\ZZ_M}(\f,\g)=\sum_{\chi \in \operatorname{Repr}(\ZZ_M)} \Omega^{fake}(\f_\chi, \g_{\chi}).\]
The Lagrangian polarization is given by $\K^{fake}_\pm\otimes \operatorname{Repr}(\ZZ_M)$, and the dilaton shift by $\vv = (1-q) v$. 

The specifics of the orbifold situation, however, is that the evaluation maps involved in the construction of the invariants take values in the {\em inertia orbifold} $IX$, in the case of the orbifold $X/\ZZ_M$ consisting of $M$ disjoint copies of $X$, which are labeled not by representations of $\ZZ_M$, but by its elements $h\in \ZZ_M$ (referred to as {\em sectors}). In sector notation
\[  \f = \sum_{\chi \in \operatorname{Repr}(\ZZ_M)} \f_\chi \chi = \sum_{h \in \ZZ_M} \f^{(h)} h,\]
where (by Fourier transform)
\[ \f^{(h)} = \sum_\chi \f_{\chi} \chi (h),\ \ \f_\chi = \frac{1}{M}
\sum_h \f^{(h)}\chi (h^{-1}).\]
Consequently,
\[ (a,b)=\frac{1}{M}\sum_h (a^{(h)},b^{(h^{-1})})^{fake},\]
the symplectic form decomposes as
\[ \Omega^{fake}_{X/\ZZ_M}(\f,\g)=\frac{1}{M}\sum_{h\in \ZZ_M}\Omega^{fake}(\f^{(h)}, \g^{(h^{-1})}),\]
the polarization spaces have the form $\oplus_{h\in \ZZ_M} \K^{fake}_{\pm} h$,
where $\K^{fake}_{+}h$ is Darboux-dual to $\K^{fake}_{-}h^{-1}$, while the dilaton shift $\vv = (1-q) \mathbf{1}$ belongs to the sector of the unit element $\mathbf{1}\in \ZZ_M$.

We will label the sectors by $M$th roots of unity $\zeta$ (primitive or not) as follows. To the element $h=h_0^{rs}$, where $h_0$ is the standard generator of $\ZZ_M$, $M=rm$, and $(s,m)=1$, we assign $\zeta (h)$ to be the {\em primitive} root of unity of order $m$ such that $\zeta^s=e^{2\pi i/m}$. Conversely, to $\zeta=e^{2\pi i t/m}$, where $m|M$, and $(t,m)=1$, we assign
$h_{(\zeta)}\in \ZZ_M$ to be $h_0^{rs}$, where $r=M/m$, and $s$ is the multiplicative inverse to $t$ modulo $m$.

\section{Formulation of the results}

We describe $\D_{X/\ZZ_M}^{tw}$ in terms of $\D_{X/\ZZ_M}^{fake}$.

The Fock space where $\D_{X/\ZZ_M}^{tw}$ lies quantizes the loop space
\[  \K^{tw}_{(M)}:=\oplus_{\zeta: \zeta^M=1} \K^{(\zeta)} \]
equipped with the symplectic form $\Omega^{tw}_{(M)}$ as follows.
Let $m=m(\zeta)$ denote the order of $\zeta$ as a {\em primitive} root of unity, and let $M=mr$. On the space $K=K^0(X)\otimes \Lambda$, introduce
a new $\Lambda$-valued pairing
\[ (a,b)^{(r)}:= \chi \left(X; a\otimes b \otimes \frac{\Eu(T_X-1)}{\Eu(\Psi^r (T_X-1))}\right).\]
Here $\Eu$ is the K-theoretic Euler class defined by $\Eu L = (1-L^{-1})=
e^{-\sum_{k>0} L^{-k}/k}$   on line bundles, and extended to arbitrary complex vector bundles by multiplicativity using the splitting principle.
The pairing satisfies
\[ (\Psi^r a, \Psi^r b)^{(r)}= r\Psi^r (a,b),\]
which is simply the abstract Grothendieck-RR formula (called also {\em Adams-RR}) for the operation $\Psi^r$ from K-theory to itself, while the factor $r$ comes from
\[ \frac{\Eu (\Psi^r1)}{\Eu (1)}=\lim_{L\to 1} \frac{1-L^{-r}}{1-L^{-1}} =r.\]

Introduce the symplectic form on $\K^{tw}_{(M)}$: 
\[ \Omega^{tw}_{(M)} (\f,\g):= \frac{1}{M}\sum_{\zeta: \ \zeta^M=1} \Res_{q=1} (\f^{(\zeta)}(q^{-1}),\g^{(\zeta^{-1})}(q))^{(r(\zeta))}\ \frac{dq}{q}.\]

To describe the polarization in $\K_{(M)}^{tw}$, introduce basis in $\K^{(\zeta)}$:
\[ \f_{k,\a}^{(\zeta)}:=\Psi^r \left(\phi^\a (q^{1/m}-1)^k\right),\ \ \g_{k,\a}^{(\zeta)}:=r \Psi^r \left( \phi_{\a}\frac{q^{k/m}}{(1-q^{1/m})^{k+1}}\right) ,\]
  where $m=m(\zeta)$, $r=r(\zeta)$, $\phi^{\a}$ runs a basis in $K^0(X)$ Poincar\'e-dual to $\phi_\a$, and $k$ run non-negative integers.
  Then $\f_{k,\a}^{\zeta}$ run a basis in the positive space of polarization, while
  $ \g_{k,\a}^{(\zeta^{-1})}$ run the Darboux-dual basis in the negative space of the polarization in question. The generating function $\D^{tw}_{X/\ZZ_M)}$ is represented by an element $\lan \D^{tw}_{X/\ZZ_M}\ran$ in the Fock space of the symplectic loop space $(\K^{tw}_{(M)}, \Omega^{tw}_{(M)})$, using this polarization, and the dilaton shift $\vv=(1-q^M)\1=\Psi^M(1-q) \1$ (in the unit sector):
  \[ \lan D^{tw}_{X/\ZZ_M}\ran ((1-q^M)\1+\t) = \D^{tw}_{X/\ZZ_M}(\t).\]
We will also assume that a quantum state does not change when the function is multiplied by a non-zero constant (so that $\lan \D \ran$ actually denotes the 1-dimensional subspace spanned by $\D$.)  
  
To state the quantum Riemann-Roch formula relating $\lan \D^{tw}_{X/\ZZ_M}\ran$ with $\lan \D^{fake}_{X/\ZZ_M}$, define operator $\square_{(M)}: \K_{(M)}^{tw}\to \K_{X/\ZZ_M}^{fake}$ acting block-diagonally by sectors:
\[ (\square_{(M)} \f)^{(\zeta)} = \square_{\zeta, r(\zeta)} (\f^{(\zeta)}), \]
where for a primitive $m$th root of unity $\eta$ and $r=1,2,3,\dots$, 
\[ \square_{\eta, r} := e^{\textstyle \sum_{k>0} \left( \frac{\Psi^{kr}(T^*_X-1)}{k(1-\eta^{-k}q^{kr/m})}-\frac{\Psi^k(T^*_X-1)}{k(1-q^k)}\right)}.\]
We claim that $\square_{(M)}$ is symplectic, i.e.
\[ \Omega_{X/\ZZ_M}^{fake}(\square_{(M)}\f, \square_{(M)} \g) = \Omega_{(M)}^{tw}(\f,\g).\]
This follows from the identity
\[ \square_{\eta,r}(q^{-1}) \square_{\eta^{-1},r}(q)=e^{\textstyle \sum_{k>0} \frac{\Psi^{kr}(T^*_X-1)-\Psi^k(T^*_X-1)}{k}} = \frac{\Eu (T_X-1)}{\Eu (\Psi^r(T_X-1))}.\]
Note that $\square_{(M)}$ respects positive spaces of our polarizations in its source and target loop spaces, but does not respect the negative ones, nor the dilaton shifts. 

\medskip

{\tt Proposition 1.}
$\lan \D_{X/\ZZ_M}^{tw} \ran = \widehat{\square}_{(M)}\, \lan \D_{X/\ZZ_M}^{fake}\ran$.

\medskip

Let us now recall the dilaton equation, which says that in the expression $\D_X^H=e^{\sum_g \h^{g-1} \F_g}$, {\em after} the dilaton shift, the functions $\F_g$ are homogeneous of degree $2-2g$ (with some anomaly for $g=1$). Namely,
\[  (\t \p_{\t} + 2 \h \p_{\h})\lan \D_X^H \ran (\t, \h) = -\frac{\eu (X)}{24} \lan \D_X^H \ran (\t, \h) .\] 
In the transition from $\lan \D_X^H$ to $\lan \D_{X/\ZZ_M}^{tw} \ran$, the homogeneity property is preserved, because our quantization formulas (from SEction 3) for quadratic Darboux monomials are homogeneous of zero degree. This allows one to recast the dependence of $\h$ (omitting the scalar factors
such as $\h^{\eu(X)/48}$) this way:
\[ \lan \D_{X/\ZZ_M}^{tw}\ran (\t ,\h, Q) = \lan \D_{X/\ZZ_M}^{tw}\ran (\frac{\t}{\sqrt{\h}}, 1, Q).\]  
  Note that $\t$ can be rewritten by sectors as $\sum_{\zeta:\ \zeta^M=1} t^{(\zeta)}h_{(\zeta)}$, where each $t_{(\zeta)}\in \K^{fake}_{+}$.
  
Now, for each {\em primitive} $m$th root of unity $\zeta$, introduce a sequence of variables $t_r^{(\zeta)}\in \K^{fake}_{+}$, where $r=1,2,3,\dots$, and define the {\em adelic} tensor product
\[ \lan \und{\D}_X \ran \, (\{ t_r^{(\zeta)} \}, \h, Q):= \bigotimes_{M=1}^{\infty} \lan \D_{X/\ZZ_M}^{tw} \ran\, (\sum_{\zeta:\ \zeta^M=1} \frac{\t_{r(\zeta)}^{(\zeta)}}{\sqrt{\h^{r(\zeta)}}} h_{(\zeta)} , 1, Q^M),\]
where for $\zeta$ of primitive order $m|M$, we put $r(\zeta)=M/m$.   

\medskip

{\tt Proposition 2.} {\em The contribution to Wick's formula for
$\lan \D_X\ran $ of the one-vertex graph (i.e. by the moduli spaces of connected quotient curves $\hat{\Sigma}$ in the notation of Section 2) is given by the logarithm $\log \lan \und{\D}_X\ran $ of adelic tensor product.}  

\medskip

The technical point in this proposition is that the dependence of the formula on $\h$ and $Q$ correctly accounts for the Euler characteristics and degrees of the {\em covering} curves $\Sigma \to \hat{\Sigma}$. 

As we have already explained in Introduction, the adelic tensor product belongs to the Fock space associated with the symplectic loop space $(\und{\K}^{\infty},\und{\Omega}^{\infty})$, which is obtained by rearranging sectors in the direct sum of the spaces $(\K_{(M)}^{tw},\Omega_{(M)}^{tw})$. 
This direct sum comes with a Lagrangian polarization inherited from those of the summands. Let us call this polarization {\em standard}.

Recall now that adelic map $\und{\ }: (\K^{\infty}, \Omega^{\infty}) \to (\und{\K}^{\infty}, \und{\Omega}^{\infty})$, defined in Introduction, is symplectic but does not respect polarizations. More precisely, the adelic image of $\K^{\infty}_{+}$ is a proper subspace in the positive space $\und{\K}_{+}^{\infty}$ of the standard polarization, while
the adelic image of $\K_{-}^{\infty}$ is Lagrangian in $\und{\K}^{\infty}$, but
does not coincide with the negative space of the standard polarization.
Let us call {\em uniform} the polarization of the adelic loop space formed by the positive space of the standard polarization and by the adelic image of $\K^{\infty}_{-}$.

\medskip

{\tt Proposition 3.} {\em The change from the standard to the uniform polarization accounts for the edges (propagators) of Wick's summation over graphs.}  

\medskip

Sections 5-9 will be dedicated to the proof of Propositions 1-3. Also, in Section 7 we will see that the adelic embedding
$\und{\ }:\K^{\infty}_{+} \subset \und{\K}^{\infty}_{+}$ of the positive spaces
of our polarizations correctly transforms the inputs $\t_r$ of $\lan \D_X\ran$  into the inputs of the adelic tensor product (they occur in the numerators of the fake holomorphic Euler characteristics in Kawasaki's RR formula). Altogether these results imply our Main Theorem:

\medskip

{\em The adelic map transforms the quantum state $\lan \und{\D}_X\ran$ into $\lan \D_X\ran$.}  

\section{Kawasaki strata}

We begin here with a detailed description of Kawasaki strata of moduli spaces of stable maps to $X$ in terms moduli spaces of stable maps to {\em orbifolds} $X/\ZZ_M$.

Let $\phi: \gS \to X$ be a stable map of a compact nodal curve (not necessarily connected) with $n$ non-singular marked points, and let $h: \gS \to \gS$ be a symmetry of this stable map (i.e. $\phi \circ h = \phi$) which is allowed to permute the marked points. Due to the stability condition, the symmetry has finite order, and therefore induces the quotient map $\hat{\phi}: \hat{\gS} \to X$ of the quotient curve $\hat{\gS}:=\gS/(h)$.
Our nearest goal is to represent the combinatorial structure of the quotient map by a certain decorated graph $\gG$. 

Let $p: \gS \to \hat{\gS}$ denote the projection of factorization.

The {\em edges} of $\gG$ correspond to {\em unbalanced nodes} of $\hat{\gS}$. For a node $\hat{\gs} \in \hat{\gS}$, denote by $r=r(\hat{\gs})$ the cardinality of its inverse image $p^{-1}(\hat{\gs})$ in $\gS$. The inverse image is an orbit of the action of $(h)$ on $\gS$, and each point $\gs$ in it
is a node of $\gS$ fixed by $h^r$. Moreover, $h^r$ preserves each of the two branches of $\gS$ at $\gs$, and acts on the tangent lines to these branches at $\gs$ by eigenvalues $\zeta_{\pm}$. The node is {\em unbalanced} if $\zeta_{+}\zeta_{-} \neq 1$.

Normalizing the quotient curve $\hat{\gS}$ at all unbalanced nodes, we obtain a collection of connected curves $\hat{\gS}_v$ which by definition correspond to {\em vertices} $v$ of graph $\gG$, and the maps $\hat{\phi}_v: \hat{\gS}_v\to X$, obtained by the restrictions of $\hat{\phi}$. Moreover, each vertex comes with the ramified $(h)$-cover $\gS_v:=p^{-1}\hat{\gS}_v \to \hat{\gS}_v$. More precisely, let $M=M_v$ be the order of $h$ on $\gS_v$. Then outside the ramification locus, $p: \gS_v \to \hat{\gS}_v$ is a principal $\ZZ_M$-bundle. This allows one to identify $\hat{\phi}_v$ with a stable map in the sense of \cite{AGV, CR, JK} to the {\em orbifold} target space $X/\ZZ_M$, the quotient of $X$ by the trivial action of the cyclic group $(h)/(h^M)$.

The moduli space of stable maps to $X/\ZZ_M$ is characterized by certain discrete invariants, which we now describe in terms of $\hat{\phi}_v$. First, it is the arithmetical {\em genus} $\hat{g}_v$ of $\hat{\gS}_v$. Next, it is the {\em degree} $\hat{d}_v$, i.e. the homology class in $H_2(X;\ZZ)$ represented by the map $\hat{\phi}_v$. Furthermore, the vertex carries marked points, which represent in $\hat{\gS}_v$ the orbits of marked points in $\gS_v$, 
ramification points which are not marked in $\gS_v$, and (the remnants in $\hat{\gS}_v$ of) the unbalanced nodes. At each such marked point $\hat{\gs}\in\hat{\gS}_v$, the {\em order} $r=r(\hat{\gs})$ of the inverse image of $\hat{\gs}$ in $\gS_v$ is defined, as well as the eigenvalue $\zeta = \zeta(\hat{\gs})$ by which the symmetry $h^r$ of $\gS_v$ acts on the tangent line at any $\gs \in p^{-1}\hat{\gs}$. Note that $\zeta$ is a primitive $m$th root of unity for some $m=M_v/r$. Therefore for some $s$ (unique $\mod m$), we have $\zeta^s = e^{2\pi i/m}$. This determines the {\em sector} of the marked point, i.e. the element, $h^{rs}$, of the cyclic group $\ZZ_{M_v}$ which acts on $T_{\gs}\gS_v$ by the generator $e^{2\pi i/m}$ of the isotropy group of $\hat{\gs}$ in the orbifold curve $\hat{\gS}$.

Thus, the Kawasaki stratum in question is characterized by the graph $\gG$ whose vertices correspond to moduli spaces of genus $\hat{g}_v$ degree $\hat{d}_v$ stable maps to $X/\ZZ_{M_v}$ with certain numbers $\hat{n}_v$ of marked points. The marked points (which are usually depicted as flags sticking out of the vertices) are decorated by the sectors (or, equivalently, primitive $m$th roots of unity $\zeta$ with $m | M_v$), while the edges pair the
{\em unbalanced} flags ($\zeta_{+}\zeta_{-}\neq 1$) of {\em the same} order: $r_{+}=M_{v_{+}}/m_{+} = M_{v_{-}}/m_{-}=r_{-}$. 

Conversely, given such a decorated graph $\gG$, one can form the corresponding Kawasaki stratum by {\em gluing} stable maps to $X/\ZZ_{M_v}$ corresponding to the vertices of $\gG$ over the diagonal constraints ($\ev_{+} = \ev_{-}$) corresponding to the edges. More precisely, each stable map to $X/\ZZ_{M_v}$ comes equipped with a principle $\ZZ_{M_v}$-bundle, possibly ramified at the markings. The generators of the groups $\ZZ_{M_v}$ define a symmetry $h$ of the total map to $X$ from the union of the covers. Since the glued marked points have the same order $r$, the covers can be glued $h$-equivariantly, resulting in stable maps to $X$ (possibly disconnected), equipped with prescribed symmetries $h$ (of order equal to the least common multiple of all $M_v$).

By applying this construction to all (possibly disconnected) decorated graphs
$\gG$, one obtains all Kawasaki strata of all moduli spaces of (possibly disconnected) stable maps to $X$.

\medskip

{\tt Remarks.} (a) When a node $\hat{\gs}$ of the curve $\hat{\gS}$ is {\em balanced}, i.e. $h^r$ fixes a node $\gs\in p^{-1}(\hat{\gs})$ but acts on the branches of $\gS$ at the node by inverse primitive $m$th roots of unity ($\zeta_{+}\zeta_{-}=1$), the stable map is deformable, at least in the virtual sense, to a non-nodal curve within the same Kawasaki stratum. The local model
of $h^r$ near $\gs$ is given by
\[ xy=\epsilon, h^r(x,y)=(\zeta_{+} x, \zeta_{-}y),\]
where $\epsilon =0$ corresponds to the nodal curve. The requirement above that the nodes corresponding to the edges of the graph are unbalanced prevents such deformations and guarantees that the stratum of symmetric maps glued according to a given graph is maximal (e.g. in the sense that $1$ does not occur as an eigenvalue of the symmetry on the virtual normal bundle to the stratum).

\begin{figure}[htb]
\begin{center}
\epsfig{file=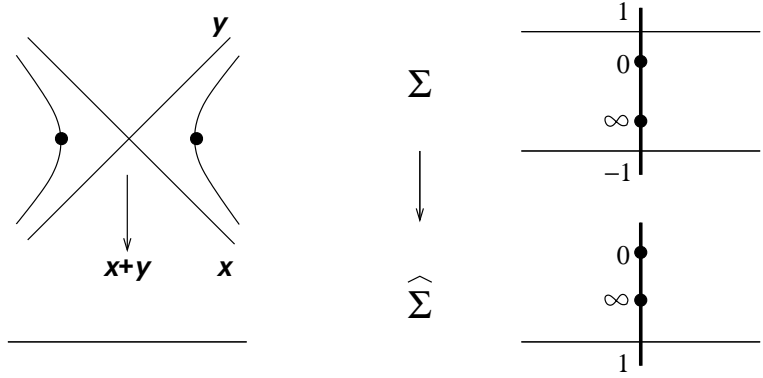} 
\end{center}
\caption{$\ZZ_2$-invariant nodes with interchanged branches}
\end{figure}

(b) One more type of deformable nodes of $\gS$ occurs when $h^r$ fixes a node
$\gs$ but interchanges the branches of $\gS$. The local model of this phenomenon can be described by the formulas:
\[ xy = \epsilon, \ h^r(x,y)=(y,x),\ \phi(x,y)=x+y,\]
so that at $\epsilon=0$, the quotient curve doesn't seem to have a node. Here is how this situation is captured in terms of orbifold stable maps. For $\epsilon \neq 0$, the map $\phi=x+y$ restricted to $xy=\epsilon$ has two ramification points: $(x,y)=\pm (\sqrt{\epsilon}, \sqrt{\epsilon})$. Thus, the
quotient curve has two marked points $\pm 2\sqrt{\epsilon}$ with inertia groups
$\ZZ_2$. When $\epsilon$ tends to $0$, the quotient curve becomes reducible,
with a new component $\CC P^1$ mapped with degree $0$, and carrying both marked points with the inertia group $\ZZ_2$ (as well as the node with the trivial isotropy group, see Figure 2). The covering curve has now $3$ components: two branches interchanged by the symmetry and connected by $\CC P^1$, which carries two marked points (say, at $z=0, \infty$), and two nodes (at $z=\pm 1$). The symmetry acts on this component by $z\mapsto -z$, so that the quotient has the node at $z^2=1$, and two marked points $z^2=0, \infty$. Thus, the quotient map, properly understood in terms of stable maps to $X/\ZZ_2$, has a balanced node of order $r=2$ with the eigenvalues $\zeta_{\pm}=1$.        

\section{Twistings}

The denominators $\str_h \wedge^{\bullet} N^*_{I\M}$ in Kawasaki's RR formula can be interpreted as certain {\em twistings} of the fake quantum K-theory of $X/\ZZ_M$, in fact a combination of several types of twistings, corresponding to different ingredients of the virtual conormal bundles.

Let $\M$ denote a Kawasaki stratum, i.e. (a component of) a moduli space $(X/\ZZ_M)_{\hat{g},\hat{n},\hat{d}}$. Let $\ft: \C \to \M$ be the corresponding universal curve, and $\ev: \C \to X/\ZZ_M$ the universal stable map, while $\tilde{\ft}: \tilde{\C}\to \M$ and $\tilde{\ev}: \tilde{\C}\to X$ denote $\ZZ_M$-equivariant lifts of $\ft$ and $\ev$ to the family of ramified $\ZZ_M$-covers.

The Kawasaki stratum $\M$ carries (the restriction to $\M$ of) the virtual tangent bundle (let's call it $\T$) to the ambient moduli space of stable maps to $X$ (say, $X_{g,n,d}$). Following \cite{Co} (see p. 99), we describe it in terms of the universal curve $\ft: \C \to \M$:
\[ \T = \tilde{\ft}_*\tilde{\ev}^*(T_X-1)\ +\ \tilde{\ft}_*(1-\tilde{L}^{-1})\ -\ (\tilde{\ft}_*\tilde{j}_*{\mathcal O}_{\tilde{\Z}})^{\vee}.\]
Here $\tilde{L}$ is the universal cotangent line bundle to the fibers of $\tilde{\ft}$ (i.e. the cotangent line bundle $L_{n+1}$ at the marked point
forgotten by $\tilde{\ft}: \C \subset X_{g,n+1,d} \to X_{g,n,d}$), and $\tilde{j}$
is the embedding of the nodal locus $\tilde{\Z}$ into $\tilde{\C}$. Loosely speaking, the three summands correspond to: (A) deformations of the maps of curves with a fixed complex structure, (B) deformations of the complex structure
of curves with fixed combinatorics, and (C) the smoothing of the nodes.

The summands carry the action of $\ZZ_M$, and can be decomposed into the eigenbundles corresponding to the eigenvalues $\la=e^{2\pi i k/M}$ of the generator. The normal bundle $N_{I\M}$, featuring in the denominator of Kawasaki's RR formula, consists of the eigenbundles corresponding to
$\la\neq 1$.

To decompose $\T$ into the eigenbundles, introduce the 1-dimensional representation $\CC_{\la}$ of $\ZZ_M$ where the generator acts by $\la$. Then the eigenbundles have the form
\begin{align*} &\T_{\la^{-1}}=(\T \otimes \CC_{\la})^{\ZZ_M}= \ft_*\ev^*[(T_X-1)\otimes\CC_{\la}]\\ &+\ft_*[(1-\tilde{L}^{-1})\otimes\ev^*\CC_{\la}]-(\ft_*[j_*{\mathcal O}_{\tilde{\Z}}\otimes \ev^* \CC_{\la}])^{\vee},\end{align*}
where $j$ is the embedding of $\Z=\tilde{\Z}/\ZZ_M$ into $\C$. 
The terms on the right are interpreted as K-theoretic push-forwards by $\ft: \C\to \M$ of orbibundles on the global quotient $\C = \tilde{\C}/\ZZ_M$. By the very definition, such push-forward automatically extracts from the sheaf cohomology its $\ZZ_M$-invariant part.

Now we use the three twisting results of \cite{ToT} to express the effect of the denominator in Kawasaki's RR formula in terms of twisted GW-invariants of orbifolds $X/\ZZ_M$.

The answer consists in the application of three operations:

(A) Transformation \[ \D_{X/\ZZ_M}^{fake} \mapsto \D_{X/\ZZ_M}^{tw}=\hat{\square}_M \D_{X/\ZZ_M}^{fake}\]
by some quantized symplectic operator (to be describe and calculated later) acting block-diagonally by $\square_M^{(h^s)}$ in the decomposition into sectors $h^s\in \ZZ_M$ of the appropriate symplectic loop spaces. 

(B) Change in the dilaton shift: $(1-q){\mathbf 1} \mapsto (1-q^M){\mathbf 1}=\Psi^M(1-q){\mathbf 1}$.

(C) Change of polarization, different on each sector (to be described later).  

In fact the three twisting theorems of \cite{ToT} are stated in terms of cohomological GW-invariants of the orbifold target ($X/\ZZ_M$ in our case). In order to relate the fake K-theory of $\M$ in Kawasaki's formula with cohomology theory, one needs to apply the three twistings with the same bundles as above, but with $\lambda=1$, and the Todd characteristic class, $\td (x)=x/(1-e^{-x})$. This results in the respective three operations described in the previous section and transforming $\D^H_{X/\ZZ_m}$ to $\D_{X/\ZZ_M}^{fake}$: by (A) application of $\hat{\qch}^{-1} \hat{\Delta}$ (the same in each sector), (B) change of the dilaton shift $-z\mapsto 1-e^z=1-q$ (in sector ${\mathbf 1}$), and (C) change of polarization from $\H_{-}$ to $\K_{-}^{fake}$ (the same in each sector). Such operations result in expressing $\D_{X/\ZZ_M}^{fake}$ in terms of $\D_X^H$ as it was explained in Section 3. The twistings A,B,C with $\lambda \neq 1$ come on the top of these, which makes it easy to phrase their outcomes directly in terms of fake quantum K-theory of $X/\ZZ_M$.

\medskip

(A) The first twisting result goes back to Tseng's ``orbifold quantum RR Theorem'' \cite{Ts}. It allows us to expresses cohomological GW-invariants of $X/\ZZ_M$ twisted by the orbibundle $E=(T_X-1)\otimes \CC_{\la}$ and by the multiplicative characteristic class $\td_{\la}$ defined by its value $1/(1-\la e^{-x})$ on a line bundle with the 1st Chern class $x$. Namely,
\[  \lan \D_{X/\ZZ_M}^{tw} \ran = \left[\prod_{k=1}^{M-1} \hat{\Delta}_{e^{2\pi i k/M}}\right] \, \lan \D_{X/\ZZ_M}^{fake} \ran,\]
where $\Delta_{e^{2\pi i k/M}}$ is the operator $\K_{X/\ZZ_M}^{fake} \to \K_{X/\ZZ_M}^{tw}$ which on the copy of $\K_X^{fake}$ corresponding to the {\em sector} $h^s\in \ZZ_M$ acts as the multiplication by the Euler--Maclaurin asymptotics of the following infinite product:
\[ \Delta_{e^{2\pi i k/M}} \sim \frac{ \prod_{l=1}^{\infty} (1-e^{2\pi i k/M} q^lq^{-\{ks/M\}})}{\prod_{l=1}^{\infty} \prod_{i=1}^{\dim_{\CC}X}(1-e^{2\pi i k/M} e^{-x_i}q^lq^{-\{ks/M\}})}.\]
Here $x_i$ are Chern roots of $T_X$, and $\{ ks/M\}$ denotes the fractional part of $ks/M$.

We rearrange the product $\prod_{k=1}^{M-1}\Delta_{e^{2\pi i k/M}}$. Let $r=(s,M)$ be the greatest common divisor of $s$ and $M$, so that $s=rs', M=rm$, $(s', m)=1$, and $ks/M=ks'/m$. Let $t'$ be inverse to $s'$ modulo $m$. Write $k=k't'+mu$ with $0\leq k'<m$. Then $\{ ks'/m \} = k'/m$ for any $u$. Since 
$\prod_{u=1}^{r}(1-Y e^{2\pi i u/r}) = 1-Y^r$ for any $Y$, we have
\begin{align*}
 \prod_{k=1}^{M-1} \prod_{l=1}^{\infty}(1-e^{2\pi i k/M} Y q^lq^{-\{ ks/M\}}) &=
 \prod_{k'} \prod_{l=1}^{\infty}(1-e^{2\pi i k't'/m} Y^r q^{rl} q^{-k'r/m}) \\
 &= \prod_{l=0}^{\infty} (1-Y^r q^{lr/m}\eta^{-l})/ \prod_{l=0}^{\infty}
 (1-Y q^l).
 \end{align*}
Here $\eta:=e^{2\pi i t'/m}$ satisfies $\eta^{s'}=e^{2\pi i /m}$, i.e. $\eta$ is the eigenvalue by which the symmetry $h^r$ acts on the tangent lines to the curves at the marked point of order $r$ and sector $h^{s}=h^{s'r}$. Also note that the Euler--Maclaurin asymptotics of the infinite product near $q=1$ is written as
\[ \prod_{l=0}^{\infty} (1-Y q^l) \sim e^{-\sum_{k>0} Y^k/k(1-q^k)}.\]
Using this, and the abbreviation $\sum_j e^{-kx_j} - 1 = \ch (\Psi^k (T^*_X-1))$,
we can summarize the above computation this way:
\[ \square_M^{(h^s)} = e^{\textstyle \sum_{k>0} \left( \frac{\Psi^{kr}(T^*_X-1)}{k(1-\eta^{-k}q^{kr/m})}-\frac{\Psi^k(T^*_X-1)}{k(1-q^k)}\right)}.\] 
The answer for $\square_M^{(h^s)}$ coincides with what was denoted by $\square_{\eta, r}$ in Section 4, where $\eta$ is a primitive $m$th root of unity, $M=mr$, $s=rs'$, and $\eta^{s'}=e^{2pi i/m}$.

\medskip

(B) The effect of the twisting by $\tilde{\ft}_*\left[(1-\tilde{L}^{-1})\otimes \tilde{\ev}^*\CC_{\la}\right]$ is described by Corollary 6.1 in \cite{ToT}. That paper, instead of the bundle $\tilde{L}$ on the covering universal curve $\tilde{\C}$, deals with the universal cotangent line bundle $L=L_{n+1}$ on $\C=\tilde{\C}/\ZZ_M$. To apply the result of that paper, it is important to realize that $\tilde{L} = p^*L$ where $p:\tilde{\C}\to \C$ is the projection of factorization. Indeed, $\tilde{L}$ is the canonical bundle of the covering curve twisted by the marked points. In local coordinates, it has a local section $x^{-1}dx$ near a marked point $x=0$, and $dx \w dy / d(xy)$ on the curves $xy=\epsilon$ near a node. The formulas
\[ \frac{dx^m}{x^m} = m \frac{dx}{x}, \ \text{and}\ \frac{dx^m \w dy^m}{d(x^my^m)} = m \frac{dx \w dy}{d(xy)}\]
identify $p^*L$ with $\tilde{L}$ near a ramified $m$-fold marked point and a balanced $m$-fold node respectively.
The answer, as we've already said, is the change of the dilaton shift:
$(1-q){\mathbf 1}\mapsto \Psi^M(1-q){\mathbf 1} = (1-q^M) {\mathbf 1}$.  

{\tt Remark.} The result does not depend on the character $\la$ of $\ZZ_M$. To understand why, the reader is invited to examine the details of the proof in \cite{ToT}, namely formula (4.2). The explanation is that the bundle $L$ is trivialized at the marked points and at the nodes (as the above local coordinate sections indicate). Consequently, Kawasaki's Chern character of $1-L^{-1}$ vanishes on all twisted sectors of the inertia orbifold $I\gS$ of the orbi-curve $\gS$. On the unit sector, however, all $\CC_{\la}$ coincide.

\medskip

(C) To describe the change of polarization caused by the twistings of type C, consider the expression $(1-L^{1/m}_{+}\otimes L^{1/m}_{-})^{-1}$. It comes from the  inverse to the K-theoretic Euler class $1-L^{1/m}_{+}L^{1/m}_{-}$ of the virtual normal line bundle  to the nodal stratum in $\M$ at the nodes of order $r=M/m$, assuming that $L_{\pm}$ represent the universal cotangent lines to the branches of quotient curve at the node. We expand the expression in powers of $L_{-}^{1/m}-1$:
\begin{align*} \frac{1}{1-L^{1/m}_{+}\otimes L^{1/m}_{-}}&=\frac{1}{1-L^{1/m}_{-}-L^{1/m}_{-}\otimes (L^{\1/m}_{+}-1)} \\
  &= \sum_{k\geq 0}\frac{L_{-}^{k/m}}{(1-L^{1/m}_{-})^{k+1}}\otimes (L^{1/m}_{+}-1)^k.
\end{align*}
Let $\{ \phi_{\ga} \}$ and $\{ \phi^{\ga} \}$ denote bases in $K^0(X)$ dual with respect to the K-theoretic Poincar{\'e} pairing. In the subspace $\K_{+}^{fake} h^{-s} \subset \K_{X/\ZZ_M}^{fake}$ (here $h^{-s}$ indicates the sector, and $r=(s,M)$ is assumed), we have a topological basis in $\K^{fake}_{+}$ ($k\geq 0$, $\ga=1,\dots, \dim K^0(X)$):
\[ \Psi^r\left( \phi^{\ga} (q^{1/m}-1)^k\right) = \Psi^r(\phi^{\ga}) (q^{r/m}-1)^k.\]
Then the following rational functions
\[ r\Psi^r\left( \phi_{\ga} \frac{q^{k/m}}{(1-q^{1/m})^{k+1}}\right) =
  r\Psi^r(\phi_{\ga}) \frac{q^{kr/m}}{(1-q^{r/m})^{k+1}},\]
 expanded into Laurent series near q=1, span the negative space of the polarization in question in the sector $h^{-s}$ of $\K^{fake}_{X/\ZZ_M}$.
  Moreover, the indicated vectors altogether form a Darboux basis in $\K_{X/\ZZ_M}^{tw}$ with respect to the symplectic form based on the following twisted pairing:
  \[ (a h^s, b h^t)^{(r)} = \frac{\delta_{h^sh^t, 1}}{M}\, \int_X \td (\Psi^r(T_X-1) \ch (a) \ch (b).\]

The result just described can be derived from a general theorem in \cite{ToT} (see Corollary 6.3 therein). It can be justified in a more direct way as well. Namely, in the non-orbifold situation, the effect of the nodal twisting leads, as it was found in the thesis \cite{Co} of T. Coates, to the change of polarization based (as it has just been described) on the ``inverse Euler class'' $(1-L_{+}\otimes L_{-})^{-1}$. In our situation of the target $X/\ZZ_M$, the smoothing of the nodes of order $r=1$ contributes into the virtual tangent bundle to Kawasaki's stratum $\M$ the same 1-dimensional summand, $\tilde{L}^{-1}_{+}\tilde{L}^{-1}_{-}=L^{-1/M}_{+}L^{-1/M}_{-}$, as into the virtual tangent bundle $\T$ to the ambient moduli space of stable maps to $X$. This means that Coates' computation still applies, with the only change that the ``inverse Euler class'' has the form $(1-L_{+}^{1/M}\otimes L_{-}^{1/M})^{-1}$. In the case of nodes of order $r>1$, the covering curves contain the $\ZZ_M$-orbit consisting of $r$ copies $\ZZ_m$-invariant nodes ($mr=M$), each contributing  into $\T$ a copy of $L_{+}^{-1/m}L_{-}^{-1/m}$, cyclically permuted by $\ZZ_r=\ZZ_M/\ZZ_m$. The ``inverse Euler class'' of their sum is $\Psi^r(1-L_{+}^{1/m}\otimes L_{-}^{1/m})^{-1}=(1-L_{+}^{r/m}\otimes L_{-}^{r/m})^{-1}$ due to the following fact that $\Psi^r(V)=\tr_h V^{\otimes r}$,
where $h$ acts on the tensor product by the cyclic permutation of the $r$ factors. 

\medskip

This completes the proof of Proposition 1.

\medskip

{\tt Remark.} We should revisit the phenomenon of $\ZZ_2$-invariant nodes with interchanged branches to examine their contribution to the type C twistings.
The cotangent line bundles $L_{\pm}$ to the branches at the node are identified by the $\ZZ_2$-symmetry: $L_{+}=L_{-}=: L$. Respectively the smoothing of the node contributes $L_{+}^{-1}L_{-}^{-1}=L^{-2}$ to the tangent bundle $\T$, and the corresponding Euler factor in the denominator of Kawasaki's formula is $1-L^2$.
It turns out that the interpretation of the situation in terms of maps to $X/\ZZ_2$ leads to the same contribution of the nodal locus. The line bundles  $L_{\pm}$ are now identified with the cotangent lines to the interchanged branches at the two nodes $\pm 1$ of the resolved curve (top right on Figure 2). Since the configuration of $0,\infty,1,-1$ on the exceptional $\CC P^1$ (vertical line at the top right) is standard, the tangent lines to this $\CC P^1$ at $\pm 1$ are trivialized. Consequently the smoothing deformation modes of the curve add up to $L_{+}^{-1}\oplus L_{-}^{-1}$ with the $\ZZ_2$-action interchanging the summands. Therefore the Euler factors representing the $\ZZ_2$-invariant and anti-invariant modes in the denominator of Kawasaki's formula are $1-L$ and $1+L$, and their product is $1-L^2$, i.e. the same as above. 

\section{Inputs}

We denote by $\t_r(q)=\sum_{k\in \ZZ} t_{r,k}q^k$ the inputs in the total descendant potential $\D_X$ of quantum K-theory on $X$, corresponding to the cycles of length $r=1,2,3,\dots$, and examine how they contribute to the numerators in Kawasaki's RR formula on the stratum $\M$ (still assuming that
the decorated graph of the stratum consists of one vertex). 

The numerators have the form of the {\em trace} $\tr_h$ of the tensor product of contributions which come from the marked points.

Let $M$ be the degree of the covers associated with a given vertex.
Let $L=L_i$ denote the universal cotangent line at a marked point on the {\em quotient} curve $\hat{\gS}=\gS /\ZZ_M$, $r=r_i$ the order of the marked point, $m=M/r$ the ramification index (of the $r$ copies in $\gS$) of this marked point, and $\zeta$ the primitive $m$th root of unity by which the symmetry $h^r$ acts on the {\em tangent} line to $\gS$ at each of the $r$ copies of the marked point. Omitting the index $i$ and the pull-back by the evaluation map $\ev_i$, we can express the resulting input of $\D_{X/\ZZ_M}^{tw}$ in the sector
determined by $\zeta$ this way: 
\[ \tr_h [\t_r (\zeta^{-1} L^{1/m})]^{\otimes r} = \Psi^r[ \t_r (\zeta^{-1} L^{1/m})]
=\sum_{k\in \ZZ}\Psi^r(t_{r,k})\zeta^{-k}L^{kr/m}.\]

Note the presence of the weight factors $1/\prod r_i^{l_i}$ in front of the supertraces $\str_h$ in the definition of the correlators involved in $\D_X$. In the expression of the correlators in terms of Kawasaki strata, these factors are compensated in the following way. Given a stable map $\hat{\Sigma} \to X/\ZZ_M$, a marked point $\hat{\sigma}\in \hat{\Sigma}$ with the ramification index $m=M/r$ represents stable maps $\Sigma \to X$ with a prescribed symmetry, which in particular cyclically permutes $r$ marked points $\sigma_1,\dots,\sigma_r$ over $\hat{\sigma}$. Even when the indices of the $r$ marked points are already decided (e.g. $1$ goes to $2$ etc. goes to $r$ goes to $1$), there still remain $r$ choices for deciding which of the marked points $\sigma_i\in \Sigma$ is numbered by $1$. Thus totally for each map $\hat{\Sigma} \to X/\ZZ_M$ in the Kawasaki stratum there are $\prod r_i^{l_i}$ symmetric maps $\Sigma \to X$, and this compensates the weight factor.   
 
\medskip

It is now time to realize that not all marked points of $\hat{\gS}$ come from  marked points of $\phi: \gS \to X$. Namely, in the theory of stable maps to $X/\ZZ_M$, all ramification points are declared marked, even if they are unmarked for the covering stable map to $X$. Consequently, the virtual cotangent bundle $\T^*$ which was analyzed in the previous section, and whose Euler class occurs in the denominator of Kawasaki's formula, is in fact the cotangent bundle to the ambient moduli space of stable maps to $X$ with extra marked points introduced at the ramifications. To compensate for these modes of deformation of stable maps, we thus need to multiply the numerator by the appropriate Euler class. Namely, if our marked point is such a ramification point, the correction has the form (one factor $1-\zeta^{-1}L^{1/m}$ per each of the $r$ copies of the ramification points): 
\[ \tr_h (1-\zeta^{-1}L^{1/m})^{\otimes r} = \Psi^r (1-\zeta^{-1} L^{1/m}) = 1-\zeta^{-1} L^{r/m}.\]
There is an exception: the unramified marked points ($m=1$, $r=M$) of $\hat{\gS}$ can come only from the orbits of marked points on $\gS$.
Note however, that in this case the same formula yields
\[ \Psi^M (1-L) = 1-L^M,\]
which agrees with the dilaton shift in $\D_{X/\ZZ_M}^{tw}$ in the unramified sector.

To summarize our observations, let us assume that the generating function $\D_{X/\ZZ_M}^{tw}$ is already dilaton-shifted by $1-q^M$ in the unit sector, and denote by $t_r^{(\zeta)} \in \K^{fake}_{+}$ the input of it through the sector indicated by the primitive $m$th root of unity $\zeta$, where $r=M/m$. Then the substitution
\[ t_r^{(\zeta)}(q) = \Psi^r\left[ 1-\zeta^{-1}q^{1/m}+\t_r(\zeta^{-1}q^{1/m}) \right] , \]
factors correctly into the numerators (and denominators) of Kawasaki's RR formula. In other words, the inputs $t_r^{(\zeta)}\in \K_{+}^{fake}$ of $\D_{X/\ZZ_M}^{tw}$ (dilaton-shifted by $1-q^M$ when $\zeta=1$) are obtained from the inputs $\t_r\in K[q,q^{-1}]$ of $\D_X$ {\em dilaton shifted by $(1-q)\1$
for  each $r=1,2,3,\dots$} by expanding $\Psi^r \t_r(\zeta^{-1}q^{1/m})$ into $q-1$-series.

This is what we claimed at the end of Section 4.

\section{Hurwitz' formula}

Here we determine the discrete characteristics of the {\em covering} of the map
$\phi: \gS \to X$, given the decorated graph $\gG$ of the quotient map $\hat{\phi}: \hat{\gS} \to X$.

The degree of $\phi$ is given by $d = \sum_v \hat{d}_v M_v$,
where $\hat{d}_v \in H_2(X;\ZZ)$ is the degree of the vertex $v$, and $M_v$ is the degree of the covering $\gS_v \to \hat{\gS}_v$.

Let us find the topological Euler characteristic $\eu$ of typical curves from the moduli spaces to which $\phi: \gS \to X$ belongs. The computation is similar to that in Hurwitz' genus formula. The vertex curves $\hat{\gS}_v$ with all the $\hat{n}_v$ ramification (i.e. marked or special) points removed have the Euler characteristics $2-2\hat{g}_v-\hat{n}_v$, which need to be multiplied by the degrees $M_v$ of the coverings. Gluing in the orbits of the ramification points of order $r_i$, $i=1,\dots, \hat{n}_v$, adds $r_i$ units for each respective orbit. Each edge $e$ of order $r_e$ corresponds to an $\ZZ_{r_e}$-orbit of (unbalanced) nodes. This subtracts
$\sum_e r_e$ units from $\eu (\gS)$, but the smoothing of all nodes subtracts $\sum_e r_r$ once more. We get  
\[ \eu = \sum_v M_v (2-2\hat{g}_v - \hat{n}_v) + \sum_v \sum_{i=1}^{\hat{n}_v} r_i
- 2\sum_e r_e.\]

Recall that our eventual goal is to represent the total descendant potential of $X$ by Kawasaki's RR formula as the sum over decorated graphs of the contributions of the respective Kawasaki strata. The contribution of the stratum represented by a given $\gG$ should be somehow obtained, starting from the product of the twisted fake potentials $\D_{X/\ZZ_{M_v}}^{tw}$, corresponding to the vertices of $\gG$, and then  ``marrying'' them appropriately by the edges. What we want to discuss now is how to dispose of the Planck constant variable $\h$ and Novikov's variables $Q$ in the vertex factors in order to achieve the correct overall occurrence of $\h$ and $Q$ in the total descendant potential of $X$.

Recall that contributions to $\D_X$ are weighted by the powers $\h^{-\eu/2}$, where $\eu$ is the Euler characteristic of the curve, connected or not, mapped to $X$. We have:
\[ -\frac{\eu}{2} = \sum_v M_v (\hat{g}_v-1) + \sum_v M_v\frac{\hat{n}_v}{2}
-\sum_v \sum_{i=1}^{\hat{n}_v} \frac{r_i}{2}
+\sum_e r_e.\]
The four terms of the sum lead to the following strategy.

(i) In each factor $\lan \D_{X/\ZZ_{M_v}}^{tw} \ran$, replace $\h$ with $\h^{M_v}$.

(ii) Replace the (dilaton-shifted) input $\t$ of the marked points with $\h^{M_v/2}\t$.

(iii) At each marked point of order $r$ divide the input by (another) factor $\h^{r/2}$.

(iv) Each ``marriage'' by an edge of order $r$ should be accompanied by the factor $\h^r$.

(v) Each monomial $Q^{\hat{d}_v}$ representing in $\lan \D_{X/\ZZ_{M_v}}^{tw}\ran $ the contributions of degree $\hat{d}$ orbicurves in $X/\ZZ_{M_v}$ should be replaced with $Q^{M_v\hat{d}_v}$.    

Now, the point is that {\em due to the homogeneity of $\lan \D_{X/\ZZ_M}^{tw}\ran $}, the steps (i) and (ii) of our strategy cancel each other, and so the steps (iii), (iv), and (v) suffice.

In particular, referring to (iii) and (v), together with the results of the previous section, we find that the vertex contribution into Wick's formula
can be described in terms of $\lan \D_{X/\ZZ_M}^{tw}(\t, \h, Q)\ran $ as the adelic product:
\[ \prod_{M=1}^{\infty} \lan \D_{X/\ZZ_M}^{tw} \ran \, \left( \sum_{\zeta:\, \zeta^M=1} \Psi^{r(\zeta)}\left[\frac{\t_{r(\zeta)}(\zeta^{-1}q^{1/m(\zeta)})}{\sqrt{\h}}\right] h_{(\zeta)}, 1, Q^M\right) ,\]
where $\t_r \in K[q,q^{-1}]$ are the arguments of $\lan \D_X \ran$.

\medskip

This completes the proof of Proposition 2.

\section{Propagators}

Recall that the edges of decorated graphs $\gG$ correspond to unbalanced nodes of the quotient curves $\hat{\gS}$. Such a node {\em of order $r$} represents $r$-tuples of nodes of the covering curve $\gS$ cyclically permuted by the symmetry $h$. On the two branches of the curve $\gS$ at such a node, $h^r$ acts with the eigenvalues $\eta_{\pm}$, which are primitive roots of unity of certain orders $m_{\pm}$. The node is {\em unbalanced} if $\eta_{+}\eta_{-}\neq 1$ (regardless of whether $m_{\pm}$ coincide or not). The effect of the unbalanced node on the contribution of the stratum $\M$ (determined by $\gG$) to Kawasaki's RR formula can be described as follows.

Let $L_{\pm}$ denote the cotangent lines to the two branches of the quotient curve $\hat{\gS}$ at the node, so that $L_{\pm}^{1/\m_{\pm}}$ denote such cotangent lines to the covering curves. The following expression
\[ \Psi^r\nabla_{\eta_{+},\eta_{-}} = \Psi^r \frac{\sum_{\a} \phi_\a\otimes \phi^\a}{1-\eta^{-1}_{+}L_{+}^{1/m_{+}} \otimes \eta^{-1}_{-}L_{-}^{1/m_{-}}} =
\frac{\sum_\a \Psi^r\phi_\a \otimes \Psi^r \phi^\a}{1-\eta^{-1}_{+}L_{+}^{r/m_{+}} \otimes \eta^{-1}_{-}L_{-}^{r/m_{-}}}\]
can be considered as an element of $K [[L_{+}^{r/m_{+}}-1]]  \otimes K [[ L_{-}^{r/m_{-}}-1]]$, where
$K=K^0(X)\otimes \Lambda$, and $\{ \phi_\a \}$ and $\{ \phi^\a \}$ are Poincar\'e-dual bases in $K^0(X)$. 
In this capacity, $\Psi^r \nabla_{\eta_{+},\eta_{-}}$ act as biderivations in the variables $\t_r^{(\eta_{\pm})}$ of the factors $\lan \D_{X/\ZZ_{M_{\pm}}}^{tw}\ran $ (with $M_{\pm}=rm_{\pm}$) in the adelic tensor product $\lan \und{\D}_X\ran $. With this notation,
Wick's summation over all graphs consists in the application to $\lan\und{\D}_X\ran$ (i.e. to the contribution of one-vertex graphs) of the following ``propagator'' (edge) operator: 
\[ \lan \und{\D}_X \ran \mapsto \exp \left[ \bigoplus_{r>0} \frac{r}{2} \h^r \Psi^r \left(\sum_{\eta_{+}\eta_{-}\neq 1}
\nabla_{\eta_{+},\eta_{-}} \right) \right]\ \lan \und{\D}_X\ran .\]
The summation sign $\oplus$ is to emphasize that the operator is block-diagonal, namely the sums with different values of $r$
act on different groups of variables, $\t_r$.    

The justification of this description is quite standard. The ingredient
$\sum_{\a} \phi_{\a}\otimes \phi^{\a}$ is responsible for the ``ungluing'' of
the diagonal constraint $\Delta \subset X\times X$ at the node. The denominator
$1-\eta_{+}^{-1} L^{1/m_{+}}_{+}\otimes \eta_{-}^{-1}L^{1/m_{-}}_{-}$ represents the
trace $\str$ (from the denominator of the Kawasaki-RR formula)
of the smoothing deformation of the curve at the node, which is normal to the
Kawasaki stratum of stable maps with the prescribed symmetry $h$. The Adams
operation $\Psi^r$ occurs at the nodes of order $r$ due to the general fact:
$\tr_h (V^{\otimes r}) = \Psi^r(V)$,
assuming that $h$ acts on the tensor product by the cyclic permutation of the
$r$ factors. The factor $r$ accounts for the number of $\ZZ_r$-equivariant ways of gluing the components of the covering curve $\Sigma$ over a node of order $r$ on the quotient curve $\hat{\Sigma}$. The factors $\h^r$ comes from the item (iv) in our strategy of the previous section to account for the change of the Euler characteristic of the covering curves $\Sigma$ under gluing at the $r$ nodes. The factor $1/2$ is due to the symmetry between $\eta_{+}$ and $\eta_{-}$. 

Our goal is to show that the application of the operator
\[ e^{\textstyle \bigoplus_{r>0} r\Psi^r \h \sum_{\eta_{+}\eta_{-}\neq 1} \nabla_{\eta_{+},\eta_{-}}/2} \]
to a function on $\und{\K}_{+}$, considered as a quantum state in the standard polarization on the adelic loop space $(\und{\K}^{\infty}, \und{\Omega}^{\infty})$, is equivalent to representing the same quantum state in the uniform polarization.

In traditional Darboux coordinate notation $\pp =\{ p_{\a}\}, \q=\{ q_{\a}\}$
a second order differential operator $\h \nabla/2  = (\h/2)\sum_{\a\b} s_{\a\b} \p_{q_\a}\p_{q_\b}$ quantizes the quadratic hamiltonian $(\pp, S\pp)/2$. The time-one map generated by the corresponding hamiltonian system
$\dot{\q}=S\pp, \dot{\pp}=0$ transforms the negative polarization space $\q=0$ into $\q=S\pp$. According to Stone-von Neumann' theorem, the operator $\exp \h \nabla /2$ intertwines the representations of the Heisenberg Lie algebra in the Fock spaces corresponding to these polarizations. Thus, we need to compute the operator $S$ in our situation, and check that the space $\q=S\pp$ is the adelic image of $\K_{-}^{\infty}$. In invariant terms, the operator $S: \und{\K}_{-}\to \und{\K}_{+}$ is computed by contracting the symmetric tensor $S \in \und{\K}_{+}\otimes \und{\K}_{+}$ using the symplectic pairing
$\und{\K}_{-}\otimes \und{\K}_{+} \to \Lambda$. 

Since our operator is block-diagonal, let us first do the computation
for the block $r=1$. Here we have the adelic space $\und{\K} = \oplus_{\zeta}\K^{(\zeta)}$, where each sector $\K^{(\zeta)}$ is isomorphic to
$\K^{fake}=K((q-1))$. It is equipped with the symplectic form
\[ \und{\Omega}(\f,\g):=\sum_{\zeta} \frac{1}{m(\zeta)} \Res_{q=1}(f^{(\zeta)}(q^{-1}), g^{(\zeta^{-1})}(q))\, \frac{dq}{q}.\]
The spaces $\und{\K}_{+}$ and $\und{\K}_{-}$ of the standard polarization are spanned respectively by (the superscript indicates the only non-zero component):
\[  \f_{k,\a}^{(\zeta)}=\phi^\a (q^{1/m}-1)^k \ \ \text{and} \ \ \g_{k,\a}^{(\zeta^{-1})}= \phi_\a \frac{q^{k/m}}{(1-q^{1/m})^{k+1}},\]
which form a Darboux basis as $k$,$\a$ and $\zeta$ run their ranges. Namely
$\und{\Omega}(\f_{k,\a}^{(\zeta)}, \g_{k,\a}^{(\zeta^{-1})})=-1$, and $=0$ in all the cases when the indices mismatch.

As it was discussed earlier,
\[ \nabla_{\eta,\zeta} = \frac{\sum_\a \phi_\a \otimes \phi^\a}{1-\eta^{-1} x^{1/m} \otimes \zeta^{-1} y^{1/n}}\ \in\ K [[ x^{1/m}-1 ]] \otimes
K [[ y^{1/n}-1 ]] \]
defines a biderivation on the space of functions on $\K_{+}^{(\eta)} \oplus \K_{+}^{(\zeta)}$. Here $\eta$ and $\zeta$ are primitive roots of unity of orders $m$ and $n$ respectively with $\eta \zeta \neq 1$ (and we write $x,y$ instead of $L_{\pm}$ used earlier).  Equivalently, $\nabla_{\eta,\zeta}$ can be considered as a bilinear form on $\K_{-}^{(\eta^{-1})}\oplus \K_{-}^{(\zeta^{-1})}$ (the symbol of the biderivation), or as a linear map $\nabla^{\zeta}_{\eta}: \K_{-}^{(\zeta^{-1})} \to \K_{+}^{(\eta)}$, which is what we want to compute.

In explicit form, the linear map $\nabla_{\eta}^{\zeta}$ is described by\footnote{The negative sign comes from $\und{\Omega}(\f_{k,\a}^{(\zeta)}, \g_{k,\a}^{(\zeta^{-1})})=-1$.}   
\[ \K_{-}^{(\zeta^{-1})} \ni f=\sum_\a f^{\a}(q) \phi_\a \mapsto - \Res_{y=1} \frac{\sum_\a \phi_\a f^\a(y)}{(1-\eta^{-1}\zeta^{-1}q^{1/m} y^{-1/n})} \frac{dy^{1/n}}{y^{1/n}}.\]  
Take $f=\phi_\a q^{k/n}/(1-q^{1/n})^{k+1}$, and put $x=y^{1/n}$. Then
\begin{align*} \nabla_{\eta}^{\zeta} f &= - \phi_a \Res_{x=1} \frac{x^k}{(1-x)^{k+1}} \frac{1}{(1-\eta^{-1}\zeta^{-1}q^{1/m}x^{-1})} \frac{dx}{x} \\
  &= \phi_\a \Res_{x=\eta^{-1}\zeta^{-1}q^{1/m}}\frac{x^k}{(1-x)^{k+1}} \frac{dx}{(x-\eta^{-1}\zeta^{-1}q^{1/m})} \\ &=\phi_{\a} \frac{(\eta^{-1}\zeta^{-1}q^{1/m})^k}{(1-\eta^{-1}\zeta^{-1}q^{1/m})^{k+1}}.\end{align*}
The last expression is interpreted as an element of $\K_{+}^{(\eta)}$ by expanding it as a power series in $q^{1/m}-1$.

Note that when $\eta$ runs all roots of unity, $k$ runs all non-negative integers, and $\phi_\a$ runs a basis of $K^0(X)$, the vector monomials $\f=\f_{\a,k,\zeta}:=\phi_\a (\zeta^{-1} q)^k/(1-\zeta^{-1} q)^{k+1}$ run a basis in $\K_{-}$. The adelic map is defined so that 
\[ \und{\f}^{(\zeta^{-1})} = \phi_\a \frac{q^{k/n}}{(1-q^{1/n})^{k+1}}, \ \ \und{\f}^{(\eta)}=
\phi_{\a} \frac{(\eta^{-1}\zeta^{-1}q^{1/m})^k}{(1-\eta^{-1}\zeta^{-1}q^{1/m})^{k+1}},\]
where $n$ and $m$ are the orders of $\zeta$ and $\eta\neq \zeta^{-1}$, and the expressions have to be expanded into Laurent series near $q=1$. The above computation shows that $\und{\f}^{(\eta)}\oplus \und{\f}^{(\zeta^{-1})}\in \K_{+}^{(\eta)}\oplus \K_{-}^{(\zeta^{-1})}$ lies in the graph of $\nabla_{\eta}^{\zeta}$. Since $\f_{\a,k,\zeta}^{(\zeta^{-1})}$ form a basis in the domain $\K_{-}^{(\zeta^{-1})}$ of $\nabla_{\eta}^{\zeta}$ when $\a$ and $k$ run their ranges, we find that $\und{\f}_{\a,k,\zeta}$ form a basis in the graph of
\[ \oplus_{\eta\neq \zeta^{-1}}\nabla_{\eta}^{\zeta}: \K_{-}^{(\zeta^{-1})} \to \und{\K}_{+},\]
and altogether form a basis in the direct sum of the graphs over $\zeta$.

\medskip

For general block $r\geq 1$, we have the adelic map: $\K^{(r)}=\K \to \und{\K}^{(r)}$, which maps $\f \in \K$ to $\Psi^r\und{\f}$. It satisfies
\[ \und{\Omega}^{(r)}(\Psi^r\und{\f}, \Psi^r\und{\g}) = \frac{\Psi^{r}}{r} \und{\Omega}(\und{\f},\und{\g})= \frac{\Psi^{r}}{r}
\Omega(\f,\g) ,\]
where $\Psi^r\Omega/r$ is the restriction of $\Omega^{\infty}$ to $\K^{(r)}$, and  $\und{\Omega}^{(r)}$ is the restriction of $\und{\Omega}^{\infty}$ to the block $\und{\K}^{(r)}$ in the total adelic space $\und{\K}^{\infty}$. It is equal to
\[ \und{\Omega}^{(r)}(\f,\g):=\frac{1}{r}\sum_{\zeta} \frac{1}{m(\zeta)} \Res_{q=1} (f^{(\zeta)}(q^{-1}),g^{(\zeta^{-1})}(q))^{(r)} \frac{dq}{q} .\]
In fact the factor $1/r$ interacts with the factor $r$ in the biderivation $r \Psi^r \nabla_{\eta,\zeta}$ in such a way that the operator from $\Psi^r(\K_{-}^{(\eta)})$ to $\und{\K}^{(r)}_{+}$  generated by it (or by the corresponding bilinear form on $\Psi^r(\K_{-}^{(\eta)}\oplus \K_{-}^{(\zeta)})$) acts as
\[ \Psi^r\left(\phi_\a \frac{q^{k/n}}{(1-q^{1/n})^{k+1}} \right) \mapsto \Psi^r \left( \phi_\a \frac{(\eta^{-1}\zeta^{-1}q^{1/m})^k}{(1-\eta^{-1}\zeta^{-1}q^{1/m})^{k+1}} \right) .\] 
Therefore the graph of the map (defined by all $r\Psi^r\nabla_{\eta,\zeta}$) from the negative space $\und{\K}_{-}^{(r)}$ of the standard polarization to $\und{\K}^{(r)}_{+}$ indeed coincides with the negative space of the uniform polarization on $\und{\K}^{(r)}$, defined as the adelic image of
$\K_{-}^{(r)}=\K_{-}$.

This completes the proof of Proposition 3, and our Main Theorem follows.


\begin{thebibliography}{10000}
  
 
\bibitem{AGV} D. Abramovich, T. Graber, A. Vistoli. {\em Algebraic orbifold quantum products.} Orbifolds in mathematics and physics (Madison, WI, 2001), pp. 124. Contemp.Math., 310. Amer. Math. Soc., Providence, RI, 2002.

  


  

\bibitem{CR} W. Chen, Y. Ruan. {\em Orbifold Gromov–Witten theory.} Orbifolds in mathematics and physics (Madison, WI, 2001), pp. 2585. Contemp. Math., 310. Amer. Math. Soc., Providence, RI, 2002.

  
\bibitem{Co} T. Coates. {\em Riemann--Roch theorems in Gromov--Witten theory.} PhD thesis, 2003, available at http://math.harvard.edu/~tomc/thesis.pdf 


\bibitem{CGL} T. Coates, A. Givental. {\em Quantum cobordisms and formal group laws.} The unity of mathematics, 155--171, Progr. Math., 244, Birkh\"auser Boston, Boston, MA, 2006.


\bibitem{Ed} D. Edidin. {\em Riemann-Roch for Deligne-Mumford stacks.}
  A celebration of algebraic geometry, pp. 241 – 266.
  Clay Math. Proc. 18 Amer. Math. Soc., Providence, RI 2013.

\bibitem{Far} C. Farsi. {\em An orbifold relative index theorem.}
J. Geom. Phys. 57 (2007), no. 8, 1653–1668.

  
  
  

  
 


  




  

\bibitem{GiTo} A. Givental, V. Tonita. {\em The Hirzebruch-Riemann--Roch theorem in true genus-0 quantum K-theory.} Preprint, arXiv:1106.3136


  

\bibitem{JK} T. Jarvis, T. Kimura. {\em Orbifold quantum cohomology of the classifying space of a finite group.} Orbifolds in mathematics and physics (Madison, WI, 2001), Contemp. Math., vol. 310, Amer. Math. Soc., Providence, RI, 2002, pp. 123--134.



\bibitem{Kaw} T. Kawasaki. {\em The Riemann-Roch theorem for complex V-manifolds.} Osaka J. Math. Volume 16, Number 1 (1979), 151-159.



\bibitem{Kr} A. Kresch. {\em On the geometry of Deligne-Mumford stacks.}
  Algebraic Geometry (Seattle, 2005), Proc. Sympos. Pure Math. 80, Part 1, Amer. Math. Soc., Providence, RI, 2009, pp. 259-271
  
  
\bibitem{YPLee} Y.-P. Lee. {\em Quantum K-theory I. Foundations.} Duke Math. J. 121 (2004), no. 3, 389-424.



\bibitem{OTT} N. O'Brien, D. Toledo, Y. L. Tong. {\em Hirzebruch-Riemann-Roch for coherent sheaves.} Amer. J. of Math., v. 103, no. 2, pp. 253--271.
  

\bibitem{ToK} V. Tonita. {\em A virtual Kawasaki Riemann–Roch formula.} Paciﬁc J. Math. 268 (2014), no. 1, 249–255. arXiv:1110.3916. 

\bibitem{ToT} V. Tonita. {\em Twisted orbifold Gromov--Witten invariants.} Nagoya Math. J. 213 (2014), 141–187, arXiv:1202.4778


\bibitem{ToH} V. Tonita. {\em A formula for the total permutation-equivariant
  K-theoretic Gromov-Witten potential.} Preprint, 13 pp., arXiv:1603.09562

\bibitem{ToTs} V. Tonita, H.-H. Tseng. {\em Quantum orbifold Hirzebruch-Riemann-Roch theorem in genus zero.} Preprint, 26 pp., arXiv:1307.0262
  
\bibitem{Ts} H.-H. Tseng. {\em Orbifold quantum Riemann-Roch, Lefschetz and Serre.} Geom. Top. 14 (2010), 1–81.  
  

\end{thebibliography}
\enddocument